\documentclass[12pt, a4paper]{article}
\usepackage{}
\usepackage{mathrsfs}
\usepackage{mathrsfs}
\usepackage{amsfonts}
\usepackage{amsfonts}
\usepackage{amsfonts}
\usepackage{mathrsfs}
\usepackage{mathrsfs}
\usepackage{mathrsfs}
\usepackage{mathrsfs}
\usepackage{mathrsfs}
\usepackage{mathrsfs}
\usepackage{mathrsfs}
\usepackage{mathrsfs}
\usepackage{mathrsfs}
\usepackage{mathrsfs}
\usepackage{mathrsfs}
\usepackage{amsfonts}
\usepackage{mathbbold}
\usepackage{amsfonts}
\usepackage{amsfonts}
\usepackage{bbm}
\usepackage{amsfonts}
\usepackage{mathrsfs}
\usepackage{mathrsfs}
\usepackage{mathrsfs}
\usepackage{mathrsfs}
\usepackage{mathrsfs}
\usepackage{mathrsfs}
\usepackage{mathrsfs}
\usepackage{mathrsfs}
\usepackage{mathrsfs}
\usepackage{mathrsfs}
\usepackage{mathrsfs}
\usepackage{mathrsfs}
\usepackage{mathrsfs}
\usepackage{mathrsfs}
\usepackage{mathrsfs}
\usepackage{amsfonts}
\usepackage{mathbbold}
\usepackage{mathrsfs}

\usepackage{latexsym}
\usepackage{graphicx}
\usepackage{amssymb}

\newtheorem{theorem}{Theorem}[section]
\newtheorem{lemma}[theorem]{Lemma}

\title{{\Large \bf  Maxima of the $Q$-index for outer-planar graphs\thanks{Supported by NSFC
(Nos. 11271315, 11171290, 11101057, 11201417).}}}
\author{Guanglong Yu$^a$\thanks{E-mail addresses:
yglong01@163.com.}
~ Shu-Guang Guo$^{a}$ ~Yarong Wu$^{b}$  ~
\\ ~ \\
{\footnotesize $^a$Department of Mathematics, Yancheng Teachers
University,}\\ {\footnotesize  Yancheng, 224002, Jiangsu, P.R. China}\\
{\footnotesize $^b$SMU college of art and science, Shanghai maritime
University, Shanghai, 200135, P.R. China}}
\date{}

\begin{document}
\maketitle

\begin{abstract}
The $Q$-$index$ of graph $G$ is the largest eigenvalue $q(G)$ of its signless Laplacian $Q(G)$. In this paper, we prove that the graph $K_{1}\nabla P_{n-1}$ has the maximal $Q$-index among all outer-planar graphs of order $n$.


\end{abstract}
\baselineskip 18.6pt

\section{Introduction}

\ \ \ \ All graphs considered in this paper are undirected and
simple, i.e. no loops or multiple edges are allowed. Given a graph $G$, $Q(G)= D(G) + A(G)$ is called the $signless$ $Laplacian$ $matrix$ of $G$, where $D(G)= \mathrm{diag}(d_{1}, d_{2},
\ldots, d_{n})$ with $d_{i}= d_{G}(v_{i})$ being the degree of
vertex $v_{i}$ $(1\leq i\leq n)$, and $A(G)$ is the adjacency matrix of $G$. The $Q$-$index$ of $G$ is the largest eigenvalue $q(G)$ of its signless Laplacian $Q(G)$. From spectral graph theory, we know that if graph $G$ is connected, there is a unit positive eigenvector (called Perron eigenvector) of $Q(G)$ corresponding to $q(G)$. In recent years, the study of the $Q$-index of a graph
attracted much attention, and the reader may consult \cite{D.P.S}-\cite{DCSKS}. On the study of the $Q$-index, a hot topic is that given a class of graphs with fixed order,
what is the maxima of the $Q$-index.

Denote by $K_{n}$, $C_{n}$, $P_{n}$ a complete graph, a cycle and a path of order $n$ respectively.
The join $G\nabla H$ of disjoint graphs $G$ and $H$ is the graph
obtained from $G\cup H$ by joining each vertex of $G$ to each vertex of $H$.
It has been conjectured in \cite{DCPR} that the graph $K_{1}\nabla P_{n-1}$ has the
maximal adjacency spectral index among all outer-planar graphs. In \cite{FYIS}, for a connected outer-planar graph of order $n\geq 2$, it has been shown that $q(G)\leq n+2$. By comparisons from some examples shown in \cite{FYIS}, it appears plausible that the graph $K_{1}\nabla P_{n-1}$ has the maximal $Q$-index among all outer-planar graphs of order $n$. In this paper, we confirm that the graph $K_{1}\nabla P_{n-1}$ really has the maximal $Q$-index among all outer-planar graphs of order $n$.

\section{Preliminary}

\ \ \ \ The reader is referred to \cite{BM, FH} for the
facts about outer-planar graphs. A graph $G$ is $outer$-$planar$ if it has a planar embedding, called $standard$ $embedding$, in which all vertices lie on the boundary of its outer face. A simple outer-planar graph is (edge) $maximal$ if no edge can be added to the graph without violating outer-planarity. In the standard embedding of a maximal outer-planar graph $G$ of order $n\geq 3$, the boundary of the outer face is a Hamiltonian
cycle (a cycle contains all vertices) of $G$, and each of the other faces is triangle. Obviously, a maximal outer-planar graph is 2-connected, and in a maximal outer-planar graph, the least vertex degree is at least $2$ (in fact, a maximal outer-planar graph has at least 2 vertices with degree $2$). From a nonmaximal outer-planar graph $G$, by adding edges to $G$, a maximal outer-planar graph $G^{'}$ can be obtained. Denote by $m(G)$ the edge number of a graph $G$. For an outer-planar graph $G$, we have $m(G)\leq 2n-3$ with equality if and only if it is maximal. From spectral graph theory, for a graph $G$, it is known that $q(G+e)> q(G)$ if $e\notin E(G)$. Consequently, when we consider the maxima of the $Q$-index among outer-planar graphs, it is sufficient to consider the maximal outer-planar graphs directly.

We introduce some notations. Denote by $V(G)$ the vertex set and $E(G)$ the edge set
for a graph $G$. If there is no ambiguity, we use $d(v)$ instead of $d_{G}(v)$. We use $\Delta$ to denote the maximum vertex degree of a graph. In a graph, the notation $v_{i}\sim v_{j}$ denotes that vertex $v_{i}$ is adjacent to
$v_{j}$. Denote by $K_{s,t}$ a complete bipartite graph with one part of size $s$ and another part of size $t$. Next we introduce some working lemmas.

\begin{lemma}{\bf \cite{SHH}}\label{le2,1} 
Let $u$ be a vertex of a maximal outer-planar graph on $n\geq 2$ vertices.
Then $\displaystyle\sum_{v\sim u}d(v)\leq n+3d(u)-4$.
\end{lemma}

\begin{lemma}{\bf \cite{R.M}}\label{le2,3} 
Let $G$ be a graph. Then $\displaystyle q(G)\leq \max_{u\in V(G)}\{d(u)+\frac{1}{d(u)}\sum_{v\sim u}d(v)\}$.
\end{lemma}

\begin{lemma}{\bf \cite{DSK}}\label{le2,4} 
Let $G$ be a connected graph containing at least one edge. Then $\displaystyle q(G)\geq \Delta+1$ with equality if and only if $G\cong K_{1, n-1}$.
\end{lemma}

\section{Main results}

\begin{lemma}\label{le3,1} 
Let $G$ be a maximal outer-planar graph with order $n\geq 6$ and $\Delta(G)\leq n-4$. Then $q(G)\leq n$.
\end{lemma}

\begin{proof}
For any vertex $u\in V(G)$, then
$$d(u)+\frac{1}{d(u)}\sum_{v\sim u}d(v)\leq d(u)+\frac{n+3d(u)-4}{d(u)}\ \ \ \ \mathrm{(by\ Lemma\ \ref{le2,1})}$$
$$ =d(u)+3+\frac{n-4}{d(u)}. \ \ $$

Let $f(x)=x+3+\frac{n-4}{x}$. It can be checked that $f(x)$ is convex.
Note that $2\leq d(u)\leq n-4$. Then $$d(u)+\frac{1}{d(u)}\sum_{v\sim u}d(v)\leq\max\{5+\frac{n-4}{2}, n\}= n.$$ By Lemma \ref{le2,3}, $q(G)\leq n$. This completes the proof.
\ \ \ \ \ $\Box$
\end{proof}

Let $\mathcal {H}=K_{1}\nabla P_{n-1}$ (see Fig. 3.1). By Lemma \ref{le2,4}, we see that $q(\mathcal {H})>n$. From this, we see that among all outer-planar graphs of order $n$, the maxima of the $Q$-index is more than $n$. Combining with Lemma \ref{le3,1}, we find that among outer-planar graphs of order $n\geq 6$, the maximal degree of the graph with the maxima of the $Q$-index is more than $n-4$. Next, we consider the outer-planar graphs of order $n\geq 6$ with $\Delta=n-3,$ $n-2$ respectively.

\setlength{\unitlength}{0.6pt}
\begin{center}
\begin{picture}(545,235)
\qbezier(21,146)(21,178)(44,201)\qbezier(44,201)(68,225)(102,225)\qbezier(102,225)(135,225)(159,201)
\qbezier(159,201)(183,178)(183,146)\qbezier(183,146)(183,113)(159,90)\qbezier(159,90)(135,67)(102,67)
\qbezier(102,67)(68,67)(44,90)\qbezier(44,90)(20,113)(21,146)
\put(99,69){\circle*{4}}
\put(182,143){\circle*{4}}
\qbezier(99,69)(140,106)(182,143)
\put(171,186){\circle*{4}}
\qbezier(171,186)(135,128)(99,69)
\put(22,133){\circle*{4}}
\qbezier(99,69)(60,101)(22,133)
\put(25,174){\circle*{4}}
\qbezier(99,69)(62,122)(25,174)
\put(47,204){\circle*{4}}
\qbezier(99,69)(73,137)(47,204)
\put(115,188){\circle*{4}}
\put(132,188){\circle*{4}}
\put(96,188){\circle*{4}}
\put(158,90){\circle*{4}}
\put(46,89){\circle*{4}}
\put(96,54){$v_{0}$}
\put(29,79){$v_{1}$}
\put(1,133){$v_{2}$}
\put(3,179){$v_{3}$}
\put(27,209){$v_{4}$}
\put(175,188){$v_{n-3}$}
\put(188,145){$v_{n-2}$}
\put(164,82){$v_{n-1}$}
\put(95,24){$\mathcal {H}$}
\qbezier(383,141)(383,171)(404,192)\qbezier(404,192)(425,214)(456,214)\qbezier(456,214)(486,214)(507,192)
\qbezier(507,192)(529,171)(529,141)\qbezier(529,141)(529,110)(507,89)\qbezier(507,89)(486,67)(456,67)
\qbezier(456,68)(425,68)(404,89)\qbezier(404,89)(383,110)(383,141)
\put(386,124){\circle*{4}}
\put(451,68){\circle*{4}}
\qbezier(386,124)(418,96)(451,68)
\put(398,185){\circle*{4}}
\qbezier(451,68)(424,127)(398,185)
\put(424,207){\circle*{4}}
\qbezier(451,68)(437,138)(424,207)
\put(491,203){\circle*{4}}
\qbezier(451,68)(471,136)(491,203)
\put(526,154){\circle*{4}}
\qbezier(451,68)(488,111)(526,154)
\put(522,113){\circle*{4}}
\qbezier(451,68)(486,91)(522,113)
\qbezier(424,207)(457,205)(491,203)
\put(498,81){\circle*{4}}
\put(513,164){\circle*{4}}
\put(400,129){\circle*{4}}
\put(463,214){\circle*{4}}
\put(400,152){\circle*{4}}
\put(400,141){\circle*{4}}
\put(505,173){\circle*{4}}
\put(498,183){\circle*{4}}
\put(456,220){$v_{0}$}
\put(414,215){$v_{1}$}
\put(446,54){$v_{k}$}
\put(409,84){\circle*{4}}
\put(383,192){$v_{2}$}
\put(345,124){$v_{k-2}$}
\put(375,76){$v_{k-1}$}
\put(501,71){$v_{k+1}$}
\put(530,111){$v_{k+2}$}
\put(530,155){$v_{k+3}$}
\put(492,210){$v_{n-1}$}
\put(445,26){$G$}
\put(208,-11){Fig. 3.1. $\mathcal {H}$, $G$}
\end{picture}
\end{center}

\begin{lemma}\label{le3,2} 
Let $G$ be a maximal outer-planar graph with order $n\geq 6$ and $\Delta(G)=d_{G}(v_{k})= n-2$ (see Fig. 3.1). Then $q(G)\leq q(\mathcal {H})$.
\end{lemma}

\begin{proof}
Because $\Delta(G)= n-2$, $2\leq k \leq n-2$. By Lemma \ref{le2,4}, we know that $q(G)> n-1\geq 5$. Let $X=(x_0, x_1, \ldots, x_{n-1})^T \in R^{n}$ be the Perron eigenvector corresponding to $q(G)$, where $x_{i}$ corresponds to vertex $v_{i}$.

Note that $$q(G)x_0=2x_0+x_1+x_{n-1},\ \ \ \  \hspace{5cm} (1)$$
$$\displaystyle q(G)x_k=(n-2)x_k+x_1+x_{n-1}+\sum_{2\leq i\leq n-2, i\neq k}x_i. \ \ \hspace{1.5cm} (2)$$ (1), (2) tell us that
$$\displaystyle q(G)x_k-q(G)x_0=(n-2)x_k-2x_0+\sum_{2\leq i\leq n-2, i\neq k} x_i,$$ $$(q(G)-2)(x_k-x_0)=(n-4)x_k+\sum_{2\leq i\leq n-2, i\neq k} x_i>0.$$ It follows immediately that $x_k>x_0$.

Note that $q(G)x_1=4x_1+x_2+x_{0}+x_{n-1}+x_k$, $q(G)x_{n-1}=4x_{n-1}+x_0+x_1+x_{n-2}+x_k$. Then
$$q(G)(x_1+x_{n-1})=5(x_1+x_{n-1})+2(x_{0}+x_k)+x_2+x_{n-2}. \hspace{2cm} (3)$$
From (1) and (2), we also get that $$\displaystyle q(G)(x_k+x_0)=(n-2)x_k+2x_0+2x_1+2x_{n-1}+\sum_{2\leq i\leq n-2, i\neq k} x_i. \hspace{1cm} (4)$$
By (4)-(3), we get that $$\displaystyle q(G)(x_k+x_0)-q(G)(x_1+x_{n-1}) \ \ \hspace{6.5cm}$$$$=(n-10)x_k+3(x_{0}+x_k)-3(x_1+x_{n-1})+3(x_k-x_{0})+\sum_{3\leq i\leq n-3, i\neq k} x_i.$$
It follows that $$\displaystyle (q(G)-3)[x_k+x_0-(x_1+x_{n-1})]=(n-10)x_k+3(x_k-x_{0})+\sum_{3\leq i\leq n-3, i\neq k} x_i. \hspace{0.5cm} (5)$$
(5) tells us that if $n\geq 10$, then $x_k+x_0>x_1+x_{n-1}$.

Let $F=G-v_{1}v_{n-1}+v_{k}v_{0}$. Note the relation between the Rayleigh quotient and the largest eigenvalue of a non-negative real symmetric matrix, and note that $X^{T}Q(F)X-X^{T}Q(G)X=(x_k+x_0)^{2}-(x_1+x_{n-1})^{2}$. It follows that if $n\geq 10$, then $q(F)>X^{T}Q(F)X> X^{T}Q(G)X=q(G)$. Because $F\cong \mathcal {H}$, if $n\geq 10$, then $q(\mathcal {H})>q(G)$.

\setlength{\unitlength}{0.63pt}
\begin{center}
\begin{picture}(665,304)
\put(49,208){\circle*{4}}
\put(10,278){\circle*{4}}
\qbezier(49,208)(29,243)(10,278)
\put(93,278){\circle*{4}}
\qbezier(49,208)(71,243)(93,278)
\put(11,230){\circle*{4}}
\put(94,233){\circle*{4}}
\put(49,303){\circle*{4}}
\put(41,182){$G_{1}$}
\put(238,207){\circle*{4}}
\put(195,262){\circle*{4}}
\put(220,294){\circle*{4}}
\qbezier(238,207)(229,251)(220,294)
\put(282,226){\circle*{4}}
\put(206,223){\circle*{4}}
\qbezier(20,175.1)(20,175.3)(20,175.2)
\put(267,293){\circle*{4}}
\qbezier(238,207)(252,250)(267,293)
\put(289,262){\circle*{4}}
\qbezier(238,207)(263,235)(289,262)
\qbezier(220,294)(213,259)(206,223)
\put(230,181){$G_{2}$}
\put(434,209){\circle*{4}}
\put(389,270){\circle*{4}}
\qbezier(434,209)(411,240)(389,270)
\put(458,298){\circle*{4}}
\qbezier(434,209)(446,254)(458,298)
\put(482,268){\circle*{4}}
\qbezier(434,209)(458,239)(482,268)
\put(476,230){\circle*{4}}
\put(414,298){\circle*{4}}
\qbezier(389,270)(423,284)(458,298)
\put(395,230){\circle*{4}}
\put(426,182){$G_{3}$}
\put(268,-9){Fig. 3.2. $G_{1}$-$G_{3}$}
\qbezier(93,278)(51,278)(10,278)
\qbezier(4,255)(4,235)(18,221)\qbezier(18,221)(32,207)(52,207)\qbezier(52,207)(71,207)(85,221)
\qbezier(85,221)(100,235)(100,255)\qbezier(4,255)(4,274)(18,288)\qbezier(18,288)(32,303)(52,303)
\qbezier(52,303)(71,303)(85,288)\qbezier(85,288)(100,274)(100,255)
\qbezier(195,253)(195,233)(209,219)\qbezier(209,219)(223,205)(243,205)\qbezier(243,205)(262,205)(276,219)
\qbezier(276,219)(291,233)(291,253)\qbezier(195,253)(195,272)(209,286)\qbezier(209,286)(223,301)(243,301)
\qbezier(243,301)(262,301)(276,286)\qbezier(276,286)(291,272)(291,253)
\qbezier(388,256)(388,236)(402,222)\qbezier(402,222)(416,208)(436,208)\qbezier(436,208)(455,208)(469,222)
\qbezier(469,222)(484,236)(484,256)\qbezier(388,256)(388,275)(402,289)\qbezier(402,289)(416,304)(436,304)
\qbezier(436,304)(455,304)(469,289)\qbezier(469,289)(484,275)(484,256)
\qbezier(569,253)(569,233)(582,219)\qbezier(582,219)(596,206)(616,206)\qbezier(616,206)(635,206)(649,219)
\qbezier(649,219)(663,233)(663,253)\qbezier(569,253)(569,272)(582,286)\qbezier(582,286)(596,300)(616,300)
\qbezier(616,300)(635,300)(649,286)\qbezier(649,286)(663,272)(663,253)
\put(577,280){\circle*{4}}
\put(654,280){\circle*{4}}
\qbezier(577,280)(615,280)(654,280)
\put(617,300){\circle*{4}}
\put(662,247){\circle*{4}}
\put(649,220){\circle*{4}}
\put(587,216){\circle*{4}}
\put(570,247){\circle*{4}}
\qbezier(570,247)(612,264)(654,280)
\qbezier(570,247)(609,234)(649,220)
\qbezier(2,104)(2,84)(15,70)\qbezier(15,70)(29,57)(49,57)\qbezier(49,57)(68,57)(82,70)
\qbezier(82,70)(96,84)(96,104)\qbezier(2,104)(2,123)(15,137)\qbezier(15,137)(29,151)(49,151)
\qbezier(49,151)(68,151)(82,137)\qbezier(82,137)(96,123)(96,104)
\put(10,132){\circle*{4}}
\put(86,132){\circle*{4}}
\qbezier(10,132)(48,132)(86,132)
\put(95,101){\circle*{4}}
\put(84,73){\circle*{4}}
\put(14,73){\circle*{4}}
\put(2,101){\circle*{4}}
\put(49,151){\circle*{4}}
\put(621,207){\circle*{4}}
\put(49,57){\circle*{4}}
\qbezier(14,73)(50,103)(86,132)
\qbezier(14,73)(54,87)(95,101)
\qbezier(14,73)(49,73)(84,73)
\qbezier(14,73)(12,103)(10,132)
\put(39,31){$G_{5}$}
\put(232,80){\circle*{4}}
\qbezier(144,103)(144,83)(157,69)\qbezier(157,69)(171,56)(191,56)\qbezier(191,56)(210,56)(224,69)
\qbezier(224,69)(238,83)(238,103)\qbezier(144,103)(144,122)(157,136)\qbezier(157,136)(171,150)(191,150)
\qbezier(191,150)(210,150)(224,136)\qbezier(224,136)(238,122)(238,103)
\put(237,103){\circle*{4}}
\put(191,57){\circle*{4}}
\put(156,71){\circle*{4}}
\put(144,101){\circle*{4}}
\put(190,149){\circle*{4}}
\put(156,135){\circle*{4}}
\put(225,135){\circle*{4}}
\qbezier(156,135)(190,135)(225,135)
\qbezier(191,57)(173,96)(156,135)
\qbezier(191,57)(208,96)(225,135)
\qbezier(191,57)(214,80)(237,103)
\put(181,30){$G_{6}$}
\qbezier(287,104)(287,84)(300,70)\qbezier(300,70)(314,57)(334,57)\qbezier(334,57)(353,57)(367,70)
\qbezier(367,70)(381,84)(381,104)\qbezier(287,104)(287,123)(300,137)\qbezier(300,137)(314,151)(334,151)
\qbezier(334,151)(353,151)(367,137)\qbezier(367,137)(381,123)(381,104)
\qbezier(191,57)(167,79)(144,101)
\put(300,137){\circle*{4}}
\put(366,137){\circle*{4}}
\qbezier(300,137)(333,137)(366,137)
\put(335,151){\circle*{4}}
\put(380,112){\circle*{4}}
\put(374,82){\circle*{4}}
\put(347,60){\circle*{4}}
\put(311,63){\circle*{4}}
\put(291,84){\circle*{4}}
\put(288,112){\circle*{4}}
\qbezier(288,112)(327,125)(366,137)
\qbezier(288,112)(334,112)(380,112)
\qbezier(288,112)(299,88)(311,63)
\qbezier(288,112)(317,86)(347,60)
\qbezier(288,112)(331,97)(374,82)
\put(324,29){$G_{7}$}
\qbezier(430,103)(430,83)(443,70)\qbezier(443,70)(456,57)(476,57)\qbezier(476,57)(495,57)(508,70)
\qbezier(508,70)(522,83)(522,103)\qbezier(430,103)(430,122)(443,135)\qbezier(443,135)(456,149)(476,149)
\qbezier(476,149)(495,149)(508,135)\qbezier(508,135)(522,122)(522,103)
\put(442,135){\circle*{4}}
\put(508,135){\circle*{4}}
\qbezier(442,135)(475,135)(508,135)
\put(475,150){\circle*{4}}
\put(520,108){\circle*{4}}
\put(515,82){\circle*{4}}
\put(500,64){\circle*{4}}
\put(472,58){\circle*{4}}
\put(441,75){\circle*{4}}
\put(430,108){\circle*{4}}
\qbezier(441,75)(441,105)(442,135)
\qbezier(441,75)(474,105)(508,135)
\qbezier(441,75)(480,92)(520,108)
\qbezier(441,75)(478,79)(515,82)
\qbezier(441,75)(470,70)(500,64)
\put(467,30){$G_{8}$}
\qbezier(575,106)(575,87)(588,74)\qbezier(588,74)(601,61)(620,61)\qbezier(620,61)(638,61)(651,74)
\qbezier(651,74)(665,87)(665,106)\qbezier(575,106)(575,124)(588,137)\qbezier(588,137)(601,151)(620,151)
\qbezier(620,151)(638,151)(651,137)\qbezier(651,137)(665,124)(665,106)
\put(610,183){$G_{4}$}
\put(587,136){\circle*{4}}
\put(652,136){\circle*{4}}
\qbezier(587,136)(619,136)(652,136)
\put(620,150){\circle*{4}}
\put(664,107){\circle*{4}}
\put(655,78){\circle*{4}}
\put(629,62){\circle*{4}}
\put(598,68){\circle*{4}}
\put(580,86){\circle*{4}}
\put(576,112){\circle*{4}}
\qbezier(598,68)(587,90)(576,112)
\qbezier(598,68)(592,102)(587,136)
\qbezier(598,68)(625,102)(652,136)
\qbezier(598,68)(631,88)(664,107)
\qbezier(598,68)(626,73)(655,78)
\put(610,30){$G_{9}$}
\qbezier(570,247)(616,247)(662,247)
\qbezier(570,247)(595,227)(621,207)
\end{picture}
\end{center}

\begin{center}
 \begin{tabular}{|c|c|c|c|c|c|c|}
\hline graph & $Q$-index & graph & $Q$-index & graph & $Q$-index \\
\hline $G_{1}$  & $6.8284$  & $G_{4}$   &  $7.9908$ & $G_{7}$ & $8.8093$    \\
\hline $G_{2}$  & $7.2571$  & $G_{5}$   &  $8.0683$ & $G_{8}$ & $8.8533$    \\
\hline $G_{3}$  & $7.3908$  & $G_{6}$   &  $8.0809$ & $G_{9}$ & $8.8611$    \\ \hline

\end{tabular}
\end{center}
\begin{center}
{Table 1. The approximation of the $Q$-index for $G_{i}$ $(1\leq i\leq 9)$}
\end{center}

It can be checked that when $n = 6$, $G\cong G_{1}$; when $n = 7$, $G\cong G_{2}$ or $G\cong G_{3}$; when when $n = 8$, $G$ is isomorphic to one in $\{G_{4}, G_{5}, G_{6}\}$; when when $n = 9$, $G$ is isomorphic to one in $\{G_{7}, G_{8}, G_{9}\}$ (see Fig. 3.2). By computation with computer, we get the approximation of the $Q$-index for each $G_{i}$ $(1\leq i\leq 9)$ (see Table 1). And by computation with computer, we get that when $n = 6$, $q(\mathcal {H})\approx 6.9576$; when $n=7$, $q(\mathcal {H})\approx 7.8099$; when $n=8$, $q(\mathcal {H})\approx 8.6925$; when $n=9$, $q(\mathcal {H})\approx 9.6007$. By a simple comparison, it follows that for each $G_{i}$ $(1\leq i\leq 9)$ of order $n$ ($6\leq n\leq 9$), $q(G_{i})<q(\mathcal {H})$.
This
completes the proof. \ \ \ \ \ $\Box$
\end{proof}

\begin{lemma}\label{le3,3} 
Let $G$ be a maximal outer-planar graph with order $n \geq 7$ and $\Delta(G) = n -3$. Then $q(G)\leq q(\mathcal {H})$.
\end{lemma}

\begin{proof}

\setlength{\unitlength}{0.63pt}
\begin{center}
\begin{picture}(651,241)
\put(91,64){\circle*{4}}
\put(33,88){\circle*{4}}
\qbezier(91,64)(62,76)(33,88)
\put(147,93){\circle*{4}}
\qbezier(91,64)(119,79)(147,93)
\put(16,159){\circle*{4}}
\qbezier(91,64)(53,112)(16,159)
\put(33,117){\circle*{4}}
\put(37,107){\circle*{4}}
\put(30,127){\circle*{4}}
\put(160,157){\circle*{4}}
\qbezier(91,64)(125,111)(160,157)
\put(149,119){\circle*{4}}
\put(152,128){\circle*{4}}
\put(145,109){\circle*{4}}
\put(32,189){\circle*{4}}
\qbezier(91,64)(61,127)(32,189)
\put(144,189){\circle*{4}}
\qbezier(91,64)(117,127)(144,189)
\qbezier(32,189)(88,189)(144,189)
\put(77,213){\circle*{4}}
\qbezier(144,189)(110,201)(77,213)
\put(58,70){\circle*{4}}
\put(122,72){\circle*{4}}
\put(314,64){\circle*{4}}
\put(255,93){\circle*{4}}
\put(375,94){\circle*{4}}
\put(243,158){\circle*{4}}
\put(260,114){\circle*{4}}
\put(263,104){\circle*{4}}
\put(256,124){\circle*{4}}
\put(386,160){\circle*{4}}
\put(376,120){\circle*{4}}
\put(379,130){\circle*{4}}
\put(372,109){\circle*{4}}
\put(263,193){\circle*{4}}
\put(372,186){\circle*{4}}
\put(297,211){\circle*{4}}
\put(281,72){\circle*{4}}
\put(348,72){\circle*{4}}
\qbezier(314,64)(284,79)(255,93)
\qbezier(314,64)(344,79)(375,94)
\qbezier(314,64)(278,111)(243,158)
\qbezier(314,64)(350,112)(386,160)
\qbezier(314,64)(343,125)(372,186)
\qbezier(372,186)(334,199)(297,211)
\put(118,209){\circle*{4}}
\put(340,209){\circle*{4}}
\put(543,202){\circle*{4}}
\qbezier(243,158)(287,172)(297,211)
\qbezier(314,64)(305,138)(297,211)
\put(571,64){\circle*{4}}
\put(626,89){\circle*{4}}
\put(492,145){\circle*{4}}
\put(513,110){\circle*{4}}
\put(517,100){\circle*{4}}
\put(509,120){\circle*{4}}
\put(647,151){\circle*{4}}
\put(630,112){\circle*{4}}
\put(634,122){\circle*{4}}
\put(626,103){\circle*{4}}
\put(497,171){\circle*{4}}
\put(639,176){\circle*{4}}
\put(511,192){\circle*{4}}
\put(540,69){\circle*{4}}
\put(604,72){\circle*{4}}
\put(626,194){\circle*{4}}
\qbezier(571,64)(531,105)(492,145)
\qbezier(571,64)(609,108)(647,151)
\qbezier(571,64)(605,120)(639,176)
\qbezier(492,145)(520,158)(511,192)
\qbezier(571,64)(541,128)(511,192)
\put(511,89){\circle*{4}}
\qbezier(571,64)(543,91)(511,89)
\qbezier(571,64)(602,93)(626,89)
\put(600,213){\circle*{4}}
\qbezier(639,176)(607,183)(600,213)
\qbezier(600,213)(585,139)(571,64)
\put(530,209){\circle*{4}}
\qbezier(571,64)(550,137)(530,209)
\put(571,220){\circle*{4}}
\qbezier(571,64)(571,142)(571,220)
\put(553,205){\circle*{4}}
\put(564,207){\circle*{4}}
\put(81,24){$D_{1}$}
\put(302,25){$D_{2}$}
\put(562,25){$D_{3}$}
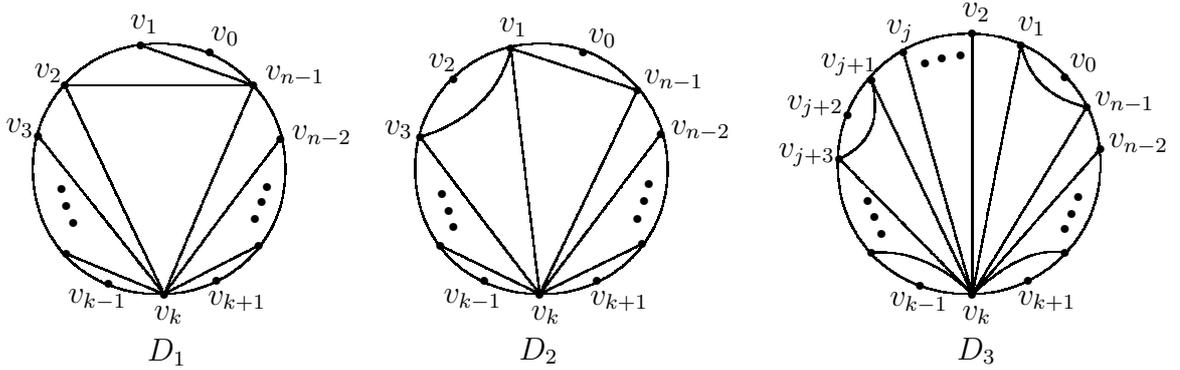
\put(235,-9){Fig. 3.3. $D_{1}$, $D_{2}$, $D_{3}$}
\put(119,216){$v_{0}$}
\put(71,223){$v_{1}$}
\put(14,194){$v_{2}$}
\put(-3,161){$v_{3}$}
\put(34,59){$v_{k-1}$}
\put(85,49){$v_{k}$}
\put(117,58){$v_{k+1}$}
\put(167,157){$v_{n-2}$}
\put(151,193){$v_{n-1}$}
\put(343,215){$v_{0}$}
\put(290,221){$v_{1}$}
\put(248,200){$v_{2}$}
\put(222,160){$v_{3}$}
\put(257,59){$v_{k-1}$}
\put(309,49){$v_{k}$}
\put(344,58){$v_{k+1}$}
\put(392,160){$v_{n-2}$}
\put(376,190){$v_{n-1}$}
\put(629,200){$v_{0}$}
\put(598,222){$v_{1}$}
\put(565,230){$v_{2}$}
\put(520,220){$v_{j}$}
\put(481,200){$v_{j+1}$}
\put(461,175){$v_{j+2}$}
\put(456,146){$v_{j+3}$}
\put(522,57){$v_{k-1}$}
\put(565,49){$v_{k}$}
\put(598,58){$v_{k+1}$}
\put(644,177){$v_{n-1}$}
\put(652,152){$v_{n-2}$}
\qbezier(491,142)(491,109)(513,86)\qbezier(513,86)(536,64)(569,64)\qbezier(569,64)(601,64)(624,86)
\qbezier(624,86)(647,109)(647,142)\qbezier(491,142)(491,174)(513,197)\qbezier(513,197)(536,220)(569,220)
\qbezier(569,220)(601,220)(624,197)\qbezier(624,197)(647,174)(647,142)
\qbezier(240,139)(240,107)(261,85)\qbezier(261,85)(283,64)(315,64)\qbezier(315,64)(346,64)(368,85)
\qbezier(368,85)(390,107)(390,139)\qbezier(240,139)(240,170)(261,192)\qbezier(261,192)(283,214)(315,214)
\qbezier(315,214)(346,214)(368,192)\qbezier(368,192)(390,170)(390,139)
\qbezier(13,139)(13,107)(34,85)\qbezier(34,85)(56,64)(88,64)\qbezier(88,64)(119,64)(141,85)
\qbezier(141,85)(163,107)(163,139)\qbezier(13,139)(13,170)(34,192)\qbezier(34,192)(56,214)(88,214)
\qbezier(88,214)(119,214)(141,192)\qbezier(141,192)(163,170)(163,139)
\end{picture}
\end{center}

Suppose $d_{G}(v_{k})=\Delta(G)$ in $G$. It can be seen that there are three cases for $G$, that is, $G\cong D_{1}$, $G\cong D_{2}$ or $G\cong D_{3}$
(see Fig. 3.3). By Lemma \ref{le2,4}, we know that $q(G) > n - 2 \geq 5$.

{\bf Case 1} $G\cong D_{1}$.

For this case, $n\geq 7$. For convenience, we suppose that $G= D_{1}$. Because $\Delta(G) = n -3$, $3\leq k\leq n-2$. Let $X=(x_0, x_1, \ldots, x_{n-1})^T \in R^{n}$ be the Perron eigenvector corresponding to $q(G)$, where $x_{i}$ corresponds to vertex $v_{i}$.

Note that $$q(G)x_0=2x_0+x_1+x_{n-1}, \ \ \hspace{3cm} (6)$$
$$q(G)x_1=3x_1+x_0+x_2+x_{n-1}.  \hspace{2.5cm} (7)$$
Then $$q(G)(x_0+x_1)=3x_0+4x_1+x_2+2x_{n-1},\hspace{3.5cm}$$  $$(q(G)-4)(x_0+x_1)=x_2+2x_{n-1}-x_0.\  \hspace{3cm} (8)$$
By (7)-(6), we get $$q(G)(x_1-x_0)=2x_1-x_0+x_2. \hspace{4.8cm}$$ Then $$(q(G)-1)(x_1-x_0)=x_1+x_{2}.\ \ \ \hspace{4cm} (9)$$
(9) implies $x_{1} > x_{0}$. Note that $$\ q(G)x_k=(n-3)x_k+\sum_{2\leq i\leq n-1, i\neq k}x_{i}. \ \ \hspace{3cm} (10)$$
By (10)+(7), we get
$$q(G)(x_k+x_1)=(n-3)x_k+x_0+3x_1+2x_2+2x_{n-1}+\sum_{3\leq i\leq n-2, i\neq k}x_{i}. \hspace{1cm} (11)$$ Then $$(q(G)-4)(x_k+x_1)=(n-7)x_k-x_1+x_0+2x_2+2x_{n-1}+\sum_{3\leq i\leq n-2, i\neq k}x_{i}.  \hspace{1cm} (12)$$ Note that $$q(G)x_2=4x_2+x_1+x_{n-1}+x_3+x_k.\ \ \ \  \hspace{3cm} (13)$$
By (13)+(6), we get that $$q(G)(x_2+x_0)=2x_0+2x_1+4x_2+x_3+x_k+2x_{n-1}.\ \ \hspace{2cm} (14)$$
By (14)-(7), we get that $$q(G)(x_2+x_0-x_1)=x_0+3x_2+x_3+x_k+x_{n-1}-x_1.$$
Then $$(q(G)-1)(x_2+x_0-x_1)=2x_2+x_3+x_k+x_{n-1}>0. \hspace{2cm} (15)$$
(15) implies that $x_2+x_0>x_1$. By (12)-(8), we get that
$$(q(G)-4)(x_k-x_0)=(n-7)x_k+2x_0+x_2-x_1+\sum_{3\leq i\leq n-2, i\neq k}x_{i}>0. \hspace{1cm} (16)$$
(16) implies $x_k>x_0$. By (10)-(7), we get that
$$q(G)(x_k-x_1)=(n-3)x_k-x_0-3x_1+\sum_{3\leq i\leq n-2, i\neq k}x_{i}.$$
Then $$(q(G)-3)(x_k-x_1)=(n-7)x_k+x_k-x_0+\sum_{3\leq i\leq n-2, i\neq k}x_{i}>0. \hspace{1.5cm} (17)$$
(17) implies that $x_k>x_1$. By (10)-(13), we get that
$$q(G)(x_k-x_2)=(n-4)x_k-3x_2-x_1+\sum_{4\leq i\leq n-2, i\neq k}x_{i}.$$
Then $$(q(G)-3)(x_k-x_2)=(n-7)x_k-x_1+\sum_{4\leq i\leq n-2, i\neq k}x_{i}. \hspace{1cm} (18)$$
(18) implies that if $n\geq 8$, then $x_k>x_2$. Note that
$$q(G)x_{n-1}=5x_{n-1}+x_0+x_1+x_2+x_k+x_{n-2}.\ \ \hspace{2.2cm} (19)$$
By (10)-(19), we get that
$$q(G)(x_k-x_{n-1})=(n-4)x_k-4x_{n-1}-x_1-x_0+\sum_{3\leq i\leq n-3, i\neq k}x_{i}.$$
Then $$(q(G)-4)(x_k-x_{n-1})=(n-8)x_k-x_1-x_0+\sum_{3\leq i\leq n-3, i\neq k}x_{i}. \hspace{1cm} (20)$$
(20) implies that if $n\geq 10$, then $x_k>x_{n-1}$. By (13)+(19), we get that
$$q(G)(x_2+x_{n-1})=x_0+2x_1+5x_2+x_3+2x_k+x_{n-2}+6x_{n-1}. \hspace{1cm} (21)$$
By (11)-(21), we get that $$q(G)(x_k+x_{1})-q(G)(x_2+x_{n-1}) \hspace{7cm}$$
$$=(n-5)x_k-3x_2-4x_{n-1}+x_{1}+\sum_{4\leq i\leq n-3, i\neq k}x_{i} \ \ \hspace{4cm}$$
$$=(n-12)x_k+3x_k-3x_2+4x_k-4x_{n-1}+x_{1}+\sum_{4\leq i\leq n-3, i\neq k}x_{i}. \hspace{1cm} (22)$$
(22) implies that if $n\geq 12$, then $x_k+x_{1}>x_2+x_{n-1}$.

Let $F = G-v_{2}v_{n-1}+v_{k}v_{1}$. Note that $X^{T}Q(F)X-X^{T}Q(G)X = (x_{k}+x_{1})^{2}-(x_{2}+x_{n-1})^{2}$.
It follows that if $n\geq 12$, then $q(F) > X^{T}Q(F)X > X^{T}Q(G)X = q(G)$. By Lemma \ref{le3,2}, it
follows immediately that if $n\geq 12$, then $q(\mathcal {H}) > q(F) > q(G)$.

When $n=7$, $G\cong G_{10}$; when $n=8$, $G$ is isomorphic to one in $\{G_{11}, G_{12}, G_{13}, G_{14}\}$; when $n=9$, $G$ is isomorphic to one in $\{G_{15}, G_{16}, G_{17}, G_{18}, G_{19}\}$; when $n=10$, $G$ is isomorphic to one in $\{G_{20}, G_{21}, G_{22}, G_{23}, G_{24}, G_{25}\}$; when $n=11$, $G$ is isomorphic to one in $\{G_{26}, G_{27}, G_{28}, G_{29}, G_{30}, G_{31}, G_{32}\}$ (see Fig. 3.4). By computation with computer, we get the approximation of the $Q$-index for each $G_{i}$ $(10\leq i\leq 32)$ (see Table 2). And by computation with computer, we get that when $n=10$, $q(\mathcal {H})\approx 10.5283$; when $n=11$, $q(\mathcal {H})\approx 11.4704$. Combining with the known results about the $Q$-index of $q(\mathcal {H})$ for order $n=7$, $8$, $9$ in the proof of Lemma \ref{le3,2}, by a simple comparison, we get that for each $G_{i}$ $(10\leq i\leq 32)$ of order $n$ ($7\leq n\leq 11$), $q(G_{i})<q(\mathcal {H})$.

\setlength{\unitlength}{0.7pt}
\begin{center}
\begin{picture}(559,644)
\put(126,625){\circle*{4}}
\put(195,625){\circle*{4}}
\put(144,639){\circle*{4}}
\put(174,639){\circle*{4}}
\put(117,594){\circle*{4}}
\put(199,583){\circle*{4}}
\put(150,555){\circle*{4}}
\qbezier(116,598)(116,579)(128,566)\qbezier(128,566)(141,554)(160,554)\qbezier(160,554)(178,554)(191,566)
\qbezier(191,566)(204,579)(204,598)\qbezier(116,598)(116,616)(128,629)\qbezier(128,629)(141,642)(160,642)
\qbezier(160,642)(178,642)(191,629)\qbezier(191,629)(204,616)(204,598)
\qbezier(126,625)(160,625)(195,625)
\qbezier(195,625)(169,632)(144,639)
\qbezier(117,594)(158,589)(199,583)
\put(246,625){\circle*{4}}
\put(315,625){\circle*{4}}
\put(264,639){\circle*{4}}
\put(294,639){\circle*{4}}
\put(237,594){\circle*{4}}
\put(322,594){\circle*{4}}
\put(247,568){\circle*{4}}
\qbezier(236,598)(236,579)(248,566)\qbezier(248,566)(261,554)(280,554)\qbezier(280,554)(298,554)(311,566)
\qbezier(311,566)(324,579)(324,598)\qbezier(236,598)(236,616)(248,629)\qbezier(248,629)(261,642)(280,642)
\qbezier(280,642)(298,642)(311,629)\qbezier(311,629)(324,616)(324,598)
\qbezier(246,625)(280,625)(315,625)
\qbezier(315,625)(289,632)(264,639)
\qbezier(237,594)(276,610)(315,625)
\qbezier(237,594)(279,594)(322,594)
\put(364,627){\circle*{4}}
\put(433,627){\circle*{4}}
\put(382,641){\circle*{4}}
\put(412,641){\circle*{4}}
\put(384,560){\circle*{4}}
\put(441,596){\circle*{4}}
\put(426,566){\circle*{4}}
\qbezier(364,627)(398,627)(433,627)
\qbezier(433,627)(407,634)(382,641)
\qbezier(384,560)(408,594)(433,627)
\qbezier(384,560)(412,578)(441,596)
\put(12,490){\circle*{4}}
\put(81,490){\circle*{4}}
\put(30,504){\circle*{4}}
\put(60,505){\circle*{4}}
\put(81,437){\circle*{4}}
\put(89,463){\circle*{4}}
\put(45,419){\circle*{4}}
\put(129,490){\circle*{4}}
\put(198,490){\circle*{4}}
\put(147,504){\circle*{4}}
\put(177,506){\circle*{4}}
\put(119,461){\circle*{4}}
\put(206,461){\circle*{4}}
\put(138,427){\circle*{4}}
\put(243,489){\circle*{4}}
\put(312,489){\circle*{4}}
\put(261,503){\circle*{4}}
\put(291,504){\circle*{4}}
\put(234,458){\circle*{4}}
\put(320,458){\circle*{4}}
\put(244,431){\circle*{4}}
\put(363,489){\circle*{4}}
\put(432,489){\circle*{4}}
\put(381,503){\circle*{4}}
\put(411,504){\circle*{4}}
\put(358,443){\circle*{4}}
\put(439,455){\circle*{4}}
\put(380,421){\circle*{4}}
\put(481,491){\circle*{4}}
\put(550,491){\circle*{4}}
\put(499,505){\circle*{4}}
\put(529,506){\circle*{4}}
\put(525,422){\circle*{4}}
\put(558,466){\circle*{4}}
\put(551,439){\circle*{4}}
\qbezier(2,463)(2,444)(14,431)\qbezier(14,431)(27,419)(46,419)\qbezier(46,419)(64,419)(77,431)
\qbezier(77,431)(90,444)(90,463)\qbezier(2,463)(2,481)(14,494)\qbezier(14,494)(27,507)(46,507)
\qbezier(46,507)(64,507)(77,494)\qbezier(77,494)(90,481)(90,463)
\qbezier(12,490)(46,490)(81,490)
\qbezier(81,490)(55,497)(30,504)
\qbezier(81,437)(81,464)(81,490)
\qbezier(129,490)(163,490)(198,490)
\qbezier(198,490)(172,497)(147,504)
\qbezier(119,461)(162,461)(206,461)
\qbezier(233,462)(233,443)(245,430)\qbezier(245,430)(258,418)(277,418)\qbezier(277,418)(295,418)(308,430)
\qbezier(308,430)(321,443)(321,462)\qbezier(233,462)(233,480)(245,493)\qbezier(245,493)(258,506)(277,506)
\qbezier(277,506)(295,506)(308,493)\qbezier(308,493)(321,480)(321,462)
\qbezier(243,489)(277,489)(312,489)
\qbezier(312,489)(286,496)(261,503)
\qbezier(234,458)(273,474)(312,489)
\qbezier(234,458)(277,458)(320,458)
\qbezier(353,462)(353,443)(365,430)\qbezier(365,430)(378,418)(397,418)\qbezier(397,418)(415,418)(428,430)
\qbezier(428,430)(441,443)(441,462)\qbezier(353,462)(353,480)(365,493)\qbezier(365,493)(378,506)(397,506)
\qbezier(397,506)(415,506)(428,493)\qbezier(428,493)(441,480)(441,462)
\qbezier(363,489)(397,489)(432,489)
\qbezier(432,489)(406,496)(381,503)
\qbezier(358,443)(395,466)(432,489)
\qbezier(358,443)(398,449)(439,455)
\qbezier(471,464)(471,445)(483,432)\qbezier(483,432)(496,420)(515,420)\qbezier(515,420)(533,420)(546,432)
\qbezier(546,432)(559,445)(559,464)\qbezier(471,464)(471,482)(483,495)\qbezier(483,495)(496,508)(515,508)
\qbezier(515,508)(533,508)(546,495)\qbezier(546,495)(559,482)(559,464)
\qbezier(481,491)(515,491)(550,491)
\qbezier(550,491)(524,498)(499,505)
\qbezier(525,422)(537,457)(550,491)
\qbezier(525,422)(541,444)(558,466)
\put(10,355){\circle*{4}}
\put(79,355){\circle*{4}}
\put(28,369){\circle*{4}}
\put(58,371){\circle*{4}}
\put(0,330){\circle*{4}}
\put(79,305){\circle*{4}}
\put(41,283){\circle*{4}}
\put(127,355){\circle*{4}}
\put(196,355){\circle*{4}}
\put(145,369){\circle*{4}}
\put(175,370){\circle*{4}}
\put(117,329){\circle*{4}}
\put(203,333){\circle*{4}}
\put(165,285){\circle*{4}}
\put(241,354){\circle*{4}}
\put(310,354){\circle*{4}}
\put(259,368){\circle*{4}}
\put(289,369){\circle*{4}}
\put(231,332){\circle*{4}}
\put(318,332){\circle*{4}}
\put(265,285){\circle*{4}}
\put(361,354){\circle*{4}}
\put(430,354){\circle*{4}}
\put(379,368){\circle*{4}}
\put(409,369){\circle*{4}}
\put(357,306){\circle*{4}}
\put(438,328){\circle*{4}}
\put(376,286){\circle*{4}}
\put(479,356){\circle*{4}}
\put(548,356){\circle*{4}}
\put(497,370){\circle*{4}}
\put(527,371){\circle*{4}}
\put(491,292){\circle*{4}}
\put(555,334){\circle*{4}}
\put(520,286){\circle*{4}}
\qbezier(0,328)(0,309)(12,296)\qbezier(12,296)(25,284)(44,284)\qbezier(44,284)(62,284)(75,296)
\qbezier(75,296)(88,309)(88,328)\qbezier(0,328)(0,346)(12,359)\qbezier(12,359)(25,372)(44,372)
\qbezier(44,372)(62,372)(75,359)\qbezier(75,359)(88,346)(88,328)
\qbezier(10,355)(44,355)(79,355)
\qbezier(79,355)(53,362)(28,369)
\qbezier(0,330)(39,318)(79,305)
\qbezier(117,328)(117,309)(129,296)\qbezier(129,296)(142,284)(161,284)\qbezier(161,284)(179,284)(192,296)
\qbezier(192,296)(205,309)(205,328)\qbezier(117,328)(117,346)(129,359)\qbezier(129,359)(142,372)(161,372)
\qbezier(161,372)(179,372)(192,359)\qbezier(192,359)(205,346)(205,328)
\qbezier(127,355)(161,355)(196,355)
\qbezier(196,355)(170,362)(145,369)
\qbezier(117,329)(160,331)(203,333)
\qbezier(231,327)(231,308)(243,295)\qbezier(243,295)(256,283)(275,283)\qbezier(275,283)(293,283)(306,295)
\qbezier(306,295)(319,308)(319,327)\qbezier(231,327)(231,345)(243,358)\qbezier(243,358)(256,371)(275,371)
\qbezier(275,371)(293,371)(306,358)\qbezier(306,358)(319,345)(319,327)
\qbezier(241,354)(275,354)(310,354)
\qbezier(310,354)(284,361)(259,368)
\qbezier(231,332)(270,343)(310,354)
\qbezier(231,332)(274,332)(318,332)
\qbezier(351,327)(351,308)(363,295)\qbezier(363,295)(376,283)(395,283)\qbezier(395,283)(413,283)(426,295)
\qbezier(426,295)(439,308)(439,327)\qbezier(351,327)(351,345)(363,358)\qbezier(363,358)(376,371)(395,371)
\qbezier(395,371)(413,371)(426,358)\qbezier(426,358)(439,345)(439,327)
\qbezier(361,354)(395,354)(430,354)
\qbezier(430,354)(404,361)(379,368)
\qbezier(357,306)(393,330)(430,354)
\qbezier(357,306)(397,317)(438,328)
\qbezier(469,329)(469,310)(481,297)\qbezier(481,297)(494,285)(513,285)\qbezier(513,285)(531,285)(544,297)
\qbezier(544,297)(557,310)(557,329)\qbezier(469,329)(469,347)(481,360)\qbezier(481,360)(494,373)(513,373)
\qbezier(513,373)(531,373)(544,360)\qbezier(544,360)(557,347)(557,329)
\qbezier(479,356)(513,356)(548,356)
\qbezier(548,356)(522,363)(497,370)
\qbezier(491,292)(519,324)(548,356)
\qbezier(491,292)(523,313)(555,334)
\put(12,222){\circle*{4}}
\put(81,222){\circle*{4}}
\put(30,236){\circle*{4}}
\put(60,237){\circle*{4}}
\put(70,159){\circle*{4}}
\put(89,202){\circle*{4}}
\put(88,181){\circle*{4}}
\put(129,222){\circle*{4}}
\put(198,222){\circle*{4}}
\put(147,236){\circle*{4}}
\put(177,237){\circle*{4}}
\put(198,171){\circle*{4}}
\put(205,197){\circle*{4}}
\put(149,154){\circle*{4}}
\put(243,221){\circle*{4}}
\put(312,221){\circle*{4}}
\put(261,235){\circle*{4}}
\put(291,236){\circle*{4}}
\put(233,200){\circle*{4}}
\put(319,200){\circle*{4}}
\put(256,156){\circle*{4}}
\put(363,221){\circle*{4}}
\put(432,221){\circle*{4}}
\put(381,235){\circle*{4}}
\put(411,236){\circle*{4}}
\put(353,201){\circle*{4}}
\put(439,201){\circle*{4}}
\put(409,152){\circle*{4}}
\put(482,224){\circle*{4}}
\put(549,223){\circle*{4}}
\put(499,237){\circle*{4}}
\put(529,238){\circle*{4}}
\put(476,178){\circle*{4}}
\put(558,203){\circle*{4}}
\put(485,162){\circle*{4}}
\qbezier(2,195)(2,176)(14,163)\qbezier(14,163)(27,151)(46,151)\qbezier(46,151)(64,151)(77,163)
\qbezier(77,163)(90,176)(90,195)\qbezier(2,195)(2,213)(14,226)\qbezier(14,226)(27,239)(46,239)
\qbezier(46,239)(64,239)(77,226)\qbezier(77,226)(90,213)(90,195)
\qbezier(12,222)(46,222)(81,222)
\qbezier(81,222)(55,229)(30,236)
\qbezier(70,159)(75,191)(81,222)
\qbezier(70,159)(79,181)(89,202)
\qbezier(119,195)(119,176)(131,163)\qbezier(131,163)(144,151)(163,151)\qbezier(163,151)(181,151)(194,163)
\qbezier(194,163)(207,176)(207,195)\qbezier(119,195)(119,213)(131,226)\qbezier(131,226)(144,239)(163,239)
\qbezier(163,239)(181,239)(194,226)\qbezier(194,226)(207,213)(207,195)
\qbezier(129,222)(163,222)(198,222)
\qbezier(198,222)(172,229)(147,236)
\qbezier(198,171)(198,197)(198,222)
\qbezier(233,194)(233,175)(245,162)\qbezier(245,162)(258,150)(277,150)\qbezier(277,150)(295,150)(308,162)
\qbezier(308,162)(321,175)(321,194)\qbezier(233,194)(233,212)(245,225)\qbezier(245,225)(258,238)(277,238)
\qbezier(277,238)(295,238)(308,225)\qbezier(308,225)(321,212)(321,194)
\qbezier(243,221)(277,221)(312,221)
\qbezier(312,221)(286,228)(261,235)
\qbezier(233,200)(276,200)(319,200)
\qbezier(353,194)(353,175)(365,162)\qbezier(365,162)(378,150)(397,150)\qbezier(397,150)(415,150)(428,162)
\qbezier(428,162)(441,175)(441,194)\qbezier(353,194)(353,212)(365,225)\qbezier(365,225)(378,238)(397,238)
\qbezier(397,238)(415,238)(428,225)\qbezier(428,225)(441,212)(441,194)
\qbezier(363,221)(397,221)(432,221)
\qbezier(432,221)(406,228)(381,235)
\qbezier(353,201)(392,211)(432,221)
\qbezier(353,201)(396,201)(439,201)
\qbezier(471,196)(471,177)(483,164)\qbezier(483,164)(496,152)(515,152)\qbezier(515,152)(533,152)(546,164)
\qbezier(546,164)(559,177)(559,196)\qbezier(471,196)(471,214)(483,227)\qbezier(483,227)(496,240)(515,240)
\qbezier(515,240)(533,240)(546,227)\qbezier(546,227)(559,214)(559,196)
\qbezier(482,224)(515,224)(549,223)
\qbezier(549,223)(524,230)(499,237)
\qbezier(476,178)(512,201)(549,223)
\qbezier(476,178)(517,191)(558,203)
\put(11,88){\circle*{4}}
\put(80,88){\circle*{4}}
\put(29,102){\circle*{4}}
\put(59,103){\circle*{4}}
\put(17,27){\circle*{4}}
\put(87,69){\circle*{4}}
\put(130,91){\circle*{4}}
\put(194,91){\circle*{4}}
\put(146,102){\circle*{4}}
\put(176,103){\circle*{4}}
\put(160,17){\circle*{4}}
\put(204,71){\circle*{4}}
\put(242,87){\circle*{4}}
\put(311,87){\circle*{4}}
\put(260,101){\circle*{4}}
\put(290,102){\circle*{4}}
\put(304,27){\circle*{4}}
\put(318,67){\circle*{4}}
\put(361,88){\circle*{4}}
\put(431,87){\circle*{4}}
\put(380,101){\circle*{4}}
\put(410,101){\circle*{4}}
\put(436,45){\circle*{4}}
\put(438,66){\circle*{4}}
\put(480,89){\circle*{4}}
\put(549,89){\circle*{4}}
\put(498,103){\circle*{4}}
\put(528,104){\circle*{4}}
\put(470,68){\circle*{4}}
\put(557,68){\circle*{4}}
\qbezier(1,61)(1,42)(13,29)\qbezier(13,29)(26,17)(45,17)\qbezier(45,17)(63,17)(76,29)\qbezier(76,29)(89,42)(89,61)
\qbezier(1,61)(1,79)(13,92)\qbezier(13,92)(26,105)(45,105)\qbezier(45,105)(63,105)(76,92)\qbezier(76,92)(89,79)(89,61)
\qbezier(11,88)(45,88)(80,88)
\qbezier(80,88)(54,95)(29,102)
\qbezier(17,27)(48,58)(80,88)
\qbezier(17,27)(52,48)(87,69)
\qbezier(118,61)(118,42)(130,29)\qbezier(130,29)(143,17)(162,17)\qbezier(162,17)(180,17)(193,29)
\qbezier(193,29)(206,42)(206,61)\qbezier(118,61)(118,79)(130,92)\qbezier(130,92)(143,105)(162,105)
\qbezier(162,105)(180,105)(193,92)\qbezier(193,92)(206,79)(206,61)
\qbezier(130,91)(162,91)(194,91)
\qbezier(194,91)(170,97)(146,102)
\qbezier(160,17)(177,54)(194,91)
\qbezier(160,17)(182,44)(204,71)
\qbezier(232,60)(232,41)(244,28)\qbezier(244,28)(257,16)(276,16)\qbezier(276,16)(294,16)(307,28)
\qbezier(307,28)(320,41)(320,60)\qbezier(232,60)(232,78)(244,91)\qbezier(244,91)(257,104)(276,104)
\qbezier(276,104)(294,104)(307,91)\qbezier(307,91)(320,78)(320,60)
\qbezier(242,87)(276,87)(311,87)
\qbezier(311,87)(285,94)(260,101)
\qbezier(304,27)(311,47)(318,67)
\qbezier(361,88)(396,88)(431,87)
\qbezier(431,87)(405,94)(380,101)
\qbezier(470,62)(470,43)(482,30)\qbezier(482,30)(495,18)(514,18)\qbezier(514,18)(532,18)(545,30)
\qbezier(545,30)(558,43)(558,62)\qbezier(470,62)(470,80)(482,93)\qbezier(482,93)(495,106)(514,106)
\qbezier(514,106)(532,106)(545,93)\qbezier(545,93)(558,80)(558,62)
\qbezier(480,89)(514,89)(549,89)
\qbezier(549,89)(523,96)(498,103)
\qbezier(470,68)(513,68)(557,68)
\qbezier(199,583)(162,604)(126,625)
\put(292,556){\circle*{4}}
\qbezier(237,594)(264,575)(292,556)
\put(356,589){\circle*{4}}
\qbezier(384,560)(374,594)(364,627)
\qbezier(354,600)(354,581)(366,568)\qbezier(366,568)(379,556)(398,556)\qbezier(398,556)(416,556)(429,568)
\qbezier(429,568)(442,581)(442,600)\qbezier(354,600)(354,618)(366,631)\qbezier(366,631)(379,644)(398,644)
\qbezier(398,644)(416,644)(429,631)\qbezier(429,631)(442,618)(442,600)
\put(5,450){\circle*{4}}
\qbezier(81,437)(43,444)(5,450)
\qbezier(81,437)(46,464)(12,490)
\put(182,423){\circle*{4}}
\qbezier(206,461)(172,444)(138,427)
\qbezier(206,461)(167,476)(129,490)
\qbezier(118,464)(118,445)(131,432)\qbezier(131,432)(144,419)(163,419)\qbezier(163,419)(181,419)(194,432)
\qbezier(194,432)(208,445)(208,464)\qbezier(118,464)(118,482)(131,495)\qbezier(131,495)(144,509)(163,509)
\qbezier(163,509)(181,509)(194,495)\qbezier(194,495)(208,482)(208,464)
\put(305,429){\circle*{4}}
\put(277,419){\circle*{4}}
\qbezier(277,419)(255,439)(234,458)
\qbezier(234,458)(269,444)(305,429)
\put(418,424){\circle*{4}}
\put(353,466){\circle*{4}}
\qbezier(358,443)(388,434)(418,424)
\qbezier(358,443)(360,466)(363,489)
\put(471,463){\circle*{4}}
\put(480,435){\circle*{4}}
\qbezier(471,463)(498,443)(525,422)
\qbezier(481,491)(503,457)(525,422)
\put(10,300){\circle*{4}}
\put(87,330){\circle*{4}}
\qbezier(79,305)(44,303)(10,300)
\qbezier(79,305)(79,330)(79,355)
\qbezier(79,305)(44,330)(10,355)
\put(127,301){\circle*{4}}
\put(196,301){\circle*{4}}
\qbezier(203,333)(184,309)(165,285)
\qbezier(203,333)(165,317)(127,301)
\qbezier(203,333)(165,344)(127,355)
\put(237,305){\circle*{4}}
\put(297,290){\circle*{4}}
\put(312,305){\circle*{4}}
\qbezier(265,285)(248,309)(231,332)
\qbezier(297,290)(264,311)(231,332)
\qbezier(312,305)(271,319)(231,332)
\put(351,328){\circle*{4}}
\put(405,286){\circle*{4}}
\put(429,301){\circle*{4}}
\qbezier(405,286)(381,296)(357,306)
\qbezier(429,301)(393,304)(357,306)
\qbezier(357,306)(359,330)(361,354)
\put(469,333){\circle*{4}}
\put(474,307){\circle*{4}}
\put(547,303){\circle*{4}}
\qbezier(491,292)(502,338)(479,356)
\qbezier(491,292)(519,298)(547,303)
\qbezier(491,292)(487,330)(469,333)
\put(34,152){\circle*{4}}
\put(2,195){\circle*{4}}
\put(13,167){\circle*{4}}
\qbezier(70,159)(41,163)(13,167)
\qbezier(70,159)(36,177)(2,195)
\qbezier(70,159)(41,191)(12,222)
\put(119,196){\circle*{4}}
\put(127,171){\circle*{4}}
\put(178,154){\circle*{4}}
\qbezier(198,171)(163,197)(129,222)
\qbezier(198,171)(158,184)(119,196)
\qbezier(198,171)(173,163)(149,154)
\qbezier(198,171)(162,171)(127,171)
\put(237,177){\circle*{4}}
\put(287,152){\circle*{4}}
\put(313,168){\circle*{4}}
\qbezier(319,200)(303,176)(287,152)
\qbezier(319,200)(287,178)(256,156)
\qbezier(319,200)(278,189)(237,177)
\qbezier(319,200)(281,211)(243,221)
\put(355,181){\circle*{4}}
\put(378,154){\circle*{4}}
\put(429,166){\circle*{4}}
\put(438,182){\circle*{4}}
\qbezier(353,201)(381,177)(409,152)
\qbezier(353,201)(391,184)(429,166)
\qbezier(353,201)(395,192)(438,182)
\put(472,203){\circle*{4}}
\put(516,152){\circle*{4}}
\put(540,161){\circle*{4}}
\put(555,179){\circle*{4}}
\qbezier(482,224)(479,201)(476,178)
\qbezier(476,178)(508,170)(540,161)
\qbezier(476,178)(515,179)(555,179)
\put(1,67){\circle*{4}}
\put(4,45){\circle*{4}}
\put(47,16){\circle*{4}}
\put(75,30){\circle*{4}}
\qbezier(17,27)(17,57)(1,67)
\qbezier(17,27)(34,71)(11,88)
\qbezier(75,30)(46,29)(17,27)
\qbezier(516,152)(496,165)(476,178)
\qbezier(353,201)(371,179)(378,154)
\put(86,47){\circle*{4}}
\qbezier(86,47)(51,37)(17,27)
\put(119,71){\circle*{4}}
\put(121,47){\circle*{4}}
\put(186,24){\circle*{4}}
\put(133,27){\circle*{4}}
\put(203,49){\circle*{4}}
\qbezier(160,17)(145,54)(130,91)
\qbezier(160,17)(139,44)(119,71)
\qbezier(160,17)(140,32)(121,47)
\qbezier(160,17)(181,33)(203,49)
\qbezier(311,87)(302,68)(304,27)
\put(283,16){\circle*{4}}
\put(317,45){\circle*{4}}
\put(235,45){\circle*{4}}
\put(250,24){\circle*{4}}
\put(233,67){\circle*{4}}
\qbezier(304,27)(268,47)(233,67)
\qbezier(304,27)(273,57)(242,87)
\qbezier(304,27)(277,26)(250,24)
\qbezier(304,27)(269,36)(235,45)
\qbezier(436,45)(414,69)(431,87)
\put(422,25){\circle*{4}}
\put(392,17){\circle*{4}}
\put(366,27){\circle*{4}}
\put(354,45){\circle*{4}}
\put(352,66){\circle*{4}}
\qbezier(436,45)(409,37)(392,17)
\qbezier(351,60)(351,41)(363,28)\qbezier(363,28)(376,16)(395,16)\qbezier(395,16)(413,16)(426,28)
\qbezier(426,28)(439,41)(439,60)\qbezier(351,60)(351,78)(363,91)\qbezier(363,91)(376,104)(395,104)
\qbezier(395,104)(413,104)(426,91)\qbezier(426,91)(439,78)(439,60)
\qbezier(436,45)(395,43)(366,27)
\qbezier(436,45)(395,45)(354,45)
\qbezier(436,45)(394,56)(352,66)
\qbezier(436,45)(398,67)(361,88)
\put(555,43){\circle*{4}}
\put(537,25){\circle*{4}}
\put(512,18){\circle*{4}}
\put(486,28){\circle*{4}}
\put(473,46){\circle*{4}}
\qbezier(557,68)(547,47)(537,25)
\qbezier(557,68)(534,43)(512,18)
\qbezier(557,68)(521,48)(486,28)
\qbezier(557,68)(515,57)(473,46)
\qbezier(557,68)(518,79)(480,89)
\put(150,530){$G_{10}$}
\put(274,530){$G_{11}$}
\put(389,530){$G_{12}$}
\put(33,397){$G_{13}$}
\put(147,397){$G_{14}$}
\put(264,397){$G_{15}$}
\put(387,397){$G_{16}$}
\put(502,397){$G_{17}$}
\put(30,258){$G_{18}$}
\put(152,258){$G_{19}$}
\put(265,258){$G_{20}$}
\put(384,258){$G_{21}$}
\put(503,258){$G_{22}$}
\put(33,129){$G_{23}$}
\put(156,129){$G_{24}$}
\put(265,129){$G_{25}$}
\put(383,129){$G_{26}$}
\put(504,129){$G_{27}$}
\put(28,-5){$G_{28}$}
\put(149,-5){$G_{29}$}
\put(264,-5){$G_{30}$}
\put(386,-5){$G_{31}$}
\put(503,-5){$G_{32}$}
\put(231,-38){Fig. 3.4 $G_{10}$-$G_{32}$}
\end{picture}
\end{center}

\bigskip

\begin{center}
 \begin{tabular}{|c|c|c|c|c|c|c|c|c|}
\hline graph & $Q$-index & graph & $Q$-index & graph & $Q$-index & graph & $Q$-index \\
\hline $G_{10}$  & $6.9895$  & $G_{16}$   &  $8.3111$ & $G_{22}$ & $9.0044$       & $G_{28}$  & $9.7983$  \\
\hline $G_{11}$  & $7.6458$  & $G_{17}$   &  $8.3225$ & $G_{23}$ & $9.0032$       & $G_{29}$     & $9.7989$ \\
\hline $G_{12}$  & $7.7873$  & $G_{18}$   &  $8.2955$ & $G_{24}$ & $8.9867$       & $G_{30}$   & $9.7977$ \\
\hline $G_{13}$  & $7.4035$  & $G_{19}$   &  $8.1101$ & $G_{25}$ & $8.8812$        & $G_{31}$   & $9.7887$  \\
\hline $G_{14}$  & $7.4641$  & $G_{20}$   &  $8.9379$ & $G_{26}$ & $9.7596$       & $G_{32}$   & $9.7274$  \\
\hline $G_{15}$  & $8.2138$  & $G_{21}$   &  $8.9954$ & $G_{27}$ & $9.7933$        &      &  \\ \hline

\end{tabular}
\end{center}
\begin{center}
{Table 2. The approximation of the $Q$-index for $G_{i}$ $(10\leq i\leq 32)$}
\end{center}

{\bf Case 2} $G\cong D_{2}$.

For this case, $n\geq 8$. For convenience, we suppose that $G= D_{2}$. Because $\Delta(G) = n -3$, $4\leq k\leq n-2$. Let $X=(x_0$, $x_1$, $\ldots$, $x_{n-1})^T \in R^{n}$ be the Perron eigenvector corresponding to $q(G)$, where $x_{i}$ corresponds to vertex $v_{i}$.

Note that $$q(G)x_0=2x_0+x_1+x_{n-1}, \hspace{4cm} (23)$$
$$q(G)x_2=2x_2+x_1+x_{3}, \ \ \hspace{4.2cm} (24)$$
$$q(G)x_k=(n-3)x_k+x_1+\sum_{3\leq i\leq n-1, i\neq k}x_{i}.\ \ \ \ \hspace{1cm} (25)$$
Then $$q(G)x_k-q(G)x_0=(n-3)x_k-2x_0+\sum_{3\leq i\leq n-2, i\neq k}x_{i},$$
$$(q(G)-2)(x_k-x_0)=(n-5)x_k+\sum_{3\leq i\leq n-2, i\neq k}x_{i}>0.$$
This implies that $x_k>x_0$.
By (25)-(24), we get that $$q(G)x_k-q(G)x_2=(n-3)x_k-2x_2+\sum_{4\leq i\leq n-1, i\neq k}x_{i}.$$
Then $$(q(G)-2)(x_k-x_2)=(n-5)x_k+\sum_{4\leq i\leq n-1, i\neq k}x_{i}>0.$$
This implies that $x_k>x_2$.
Note that $$q(G)x_1=5x_1+x_0+x_2+x_3+x_k+x_{n-1},\ \ \  \hspace{2cm} (26)$$
$$q(G)x_{n-1}=4x_{n-1}+x_0+x_1+x_{k}+x_{n-2}.\ \ \ \  \hspace{2cm} (27)$$
By (25)-(26), we get $$q(G)x_k-q(G)x_1=(n-4)x_k-4x_1-x_0-x_2+\sum_{4\leq i\leq n-2, i\neq k}x_{i}.$$
Then $$(q(G)-4)(x_k-x_1)=(n-10)x_k+2x_k-x_0-x_2+\sum_{4\leq i\leq n-2, i\neq k}x_{i}.$$
This implies that if $n\geq 10$, then $x_k>x_1$.
By (25)-(27), we get $$q(G)x_k-q(G)x_{n-1}=(n-4)x_k-3x_{n-1}-x_0+\sum_{3\leq i\leq n-3, i\neq k}x_{i}.$$
Then $$(q(G)-3)(x_k-x_{n-1})=(n-8)x_k+x_k-x_0+\sum_{3\leq i\leq n-3, i\neq k}x_{i}.$$
This implies that if $n\geq 8$, then $x_k>x_{n-1}$.
By (23)+(25), we get that
$$q(G)(x_k+x_0)=(n-3)x_k+2x_0+2x_1+2x_{n-1}+\sum_{3\leq i\leq n-2, i\neq k}x_{i}.\ \  \hspace{1cm} (28)$$
By (26)+(27), we get that
$$q(G)(x_1+x_{n-1})=2x_k+2x_0+6x_1+5x_{n-1}+x_2+x_3+x_{n-2}.\   \hspace{1.5cm} (29)$$
By (28)-(29), we get that
$$q(G)(x_k+x_0)-q(G)(x_1+x_{n-1}) \hspace{8cm}$$
$$=(n-5)x_k-4x_1-3x_{n-1}-x_2+\sum_{4\leq i\leq n-3, i\neq k}x_{i} \ \ \hspace{5cm}$$
$$=(n-13)x_k+4x_k-4x_1+3x_k-3x_{n-1}+x_k-x_2+\sum_{4\leq i\leq n-3, i\neq k}x_{i}. \hspace{1cm}(30)$$
(30) implies that if $n\geq 13$, then $x_k+x_0>x_1+x_{n-1}$.

Let $F=G-v_{1}v_{n-1}+v_{k}v_{0}$. Note that $X^{T}Q(F)X-X^{T}Q(G)X=(x_k+x_0)^{2}-(x_1+x_{n-1})^{2}$. It follows that if $n\geq 13$, then $q(F)>X^{T}Q(F)X> X^{T}Q(G)X=q(G)$. By Lemma \ref{le3,2}, it follows immediately that if $n\geq 13$, then $q(\mathcal {H})>q(F)>q(G)$.

\setlength{\unitlength}{0.7pt}
\begin{center}
\begin{picture}(558,435)
\put(73,421){\circle*{4}}
\put(146,408){\circle*{4}}
\put(104,434){\circle*{4}}
\put(132,426){\circle*{4}}
\put(64,374){\circle*{4}}
\put(146,372){\circle*{4}}
\put(62,399){\circle*{4}}
\qbezier(62,390)(62,371)(74,358)\qbezier(74,358)(87,346)(106,346)\qbezier(106,346)(124,346)(137,358)
\qbezier(137,358)(150,371)(150,390)\qbezier(62,390)(62,408)(74,421)\qbezier(74,421)(87,434)(106,434)
\qbezier(106,434)(124,434)(137,421)\qbezier(137,421)(150,408)(150,390)
\qbezier(146,408)(125,421)(104,434)
\qbezier(64,374)(105,391)(146,408)
\put(179,407){\circle*{4}}
\put(262,403){\circle*{4}}
\put(196,427){\circle*{4}}
\put(247,425){\circle*{4}}
\put(180,371){\circle*{4}}
\put(258,367){\circle*{4}}
\put(219,344){\circle*{4}}
\qbezier(176,389)(176,370)(188,357)\qbezier(188,357)(201,345)(220,345)\qbezier(220,345)(238,345)(251,357)
\qbezier(251,357)(264,370)(264,389)\qbezier(176,389)(176,407)(188,420)\qbezier(188,420)(201,433)(220,433)
\qbezier(220,433)(238,433)(251,420)\qbezier(251,420)(264,407)(264,389)
\put(296,406){\circle*{4}}
\put(379,408){\circle*{4}}
\put(310,426){\circle*{4}}
\put(364,426){\circle*{4}}
\put(295,378){\circle*{4}}
\put(381,381){\circle*{4}}
\put(319,351){\circle*{4}}
\qbezier(294,391)(294,372)(306,359)\qbezier(306,359)(319,347)(338,347)\qbezier(338,347)(356,347)(369,359)
\qbezier(369,359)(382,372)(382,391)\qbezier(294,391)(294,409)(306,422)\qbezier(306,422)(319,435)(338,435)
\qbezier(338,435)(356,435)(369,422)\qbezier(369,422)(382,409)(382,391)
\put(417,409){\circle*{4}}
\put(498,410){\circle*{4}}
\put(432,427){\circle*{4}}
\put(483,427){\circle*{4}}
\put(444,350){\circle*{4}}
\put(501,387){\circle*{4}}
\put(487,358){\circle*{4}}
\put(5,264){\circle*{4}}
\put(85,265){\circle*{4}}
\put(20,283){\circle*{4}}
\put(69,283){\circle*{4}}
\put(80,219){\circle*{4}}
\put(89,246){\circle*{4}}
\put(16,212){\circle*{4}}
\put(119,263){\circle*{4}}
\put(204,264){\circle*{4}}
\put(136,283){\circle*{4}}
\put(187,284){\circle*{4}}
\put(119,231){\circle*{4}}
\put(207,239){\circle*{4}}
\put(137,209){\circle*{4}}
\put(236,266){\circle*{4}}
\put(316,265){\circle*{4}}
\put(253,283){\circle*{4}}
\put(302,281){\circle*{4}}
\put(232,241){\circle*{4}}
\put(320,241){\circle*{4}}
\put(243,214){\circle*{4}}
\put(357,267){\circle*{4}}
\put(435,267){\circle*{4}}
\put(371,282){\circle*{4}}
\put(420,282){\circle*{4}}
\put(355,226){\circle*{4}}
\put(440,245){\circle*{4}}
\put(371,208){\circle*{4}}
\put(475,269){\circle*{4}}
\put(552,270){\circle*{4}}
\put(488,283){\circle*{4}}
\put(537,285){\circle*{4}}
\put(525,205){\circle*{4}}
\put(558,249){\circle*{4}}
\put(550,222){\circle*{4}}
\qbezier(1,246)(1,227)(13,214)\qbezier(13,214)(26,202)(45,202)\qbezier(45,202)(63,202)(76,214)
\qbezier(76,214)(89,227)(89,246)\qbezier(1,246)(1,264)(13,277)\qbezier(13,277)(26,290)(45,290)
\qbezier(45,290)(63,290)(76,277)\qbezier(76,277)(89,264)(89,246)
\qbezier(232,245)(232,226)(244,213)\qbezier(244,213)(257,201)(276,201)\qbezier(276,201)(294,201)(307,213)
\qbezier(307,213)(320,226)(320,245)\qbezier(232,245)(232,263)(244,276)\qbezier(244,276)(257,289)(276,289)
\qbezier(276,289)(294,289)(307,276)\qbezier(307,276)(320,263)(320,245)
\qbezier(352,245)(352,226)(364,213)\qbezier(364,213)(377,201)(396,201)\qbezier(396,201)(414,201)(427,213)
\qbezier(427,213)(440,226)(440,245)\qbezier(352,245)(352,263)(364,276)\qbezier(364,276)(377,289)(396,289)
\qbezier(396,289)(414,289)(427,276)\qbezier(427,276)(440,263)(440,245)
\qbezier(470,247)(470,228)(482,215)\qbezier(482,215)(495,203)(514,203)\qbezier(514,203)(532,203)(545,215)
\qbezier(545,215)(558,228)(558,247)\qbezier(470,247)(470,265)(482,278)\qbezier(482,278)(495,291)(514,291)
\qbezier(514,291)(532,291)(545,278)\qbezier(545,278)(558,265)(558,247)
\put(6,129){\circle*{4}}
\put(81,128){\circle*{4}}
\put(20,143){\circle*{4}}
\put(68,142){\circle*{4}}
\put(0,107){\circle*{4}}
\put(82,82){\circle*{4}}
\put(26,63){\circle*{4}}
\put(122,129){\circle*{4}}
\put(197,129){\circle*{4}}
\put(138,144){\circle*{4}}
\put(185,142){\circle*{4}}
\put(116,108){\circle*{4}}
\put(203,110){\circle*{4}}
\put(166,62){\circle*{4}}
\put(235,127){\circle*{4}}
\put(314,124){\circle*{4}}
\put(253,143){\circle*{4}}
\put(299,140){\circle*{4}}
\put(230,109){\circle*{4}}
\put(317,104){\circle*{4}}
\put(233,85){\circle*{4}}
\put(355,126){\circle*{4}}
\put(432,127){\circle*{4}}
\put(370,142){\circle*{4}}
\put(416,143){\circle*{4}}
\put(354,83){\circle*{4}}
\put(438,105){\circle*{4}}
\put(369,67){\circle*{4}}
\put(476,133){\circle*{4}}
\put(547,133){\circle*{4}}
\put(492,146){\circle*{4}}
\put(535,144){\circle*{4}}
\put(485,70){\circle*{4}}
\put(556,114){\circle*{4}}
\put(538,70){\circle*{4}}
\qbezier(116,105)(116,86)(128,73)\qbezier(128,73)(141,61)(160,61)\qbezier(160,61)(178,61)(191,73)
\qbezier(191,73)(204,86)(204,105)\qbezier(116,105)(116,123)(128,136)\qbezier(128,136)(141,149)(160,149)
\qbezier(160,149)(178,149)(191,136)\qbezier(191,136)(204,123)(204,105)
\qbezier(230,104)(230,85)(242,72)\qbezier(242,72)(255,60)(274,60)\qbezier(274,60)(292,60)(305,72)
\qbezier(305,72)(318,85)(318,104)\qbezier(230,104)(230,122)(242,135)\qbezier(242,135)(255,148)(274,148)
\qbezier(274,148)(292,148)(305,135)\qbezier(305,135)(318,122)(318,104)
\qbezier(349,104)(349,85)(362,72)\qbezier(362,72)(375,59)(394,59)\qbezier(394,59)(412,59)(425,72)
\qbezier(425,72)(439,85)(439,104)\qbezier(349,104)(349,122)(362,135)\qbezier(362,135)(375,149)(394,149)
\qbezier(394,149)(412,149)(425,135)\qbezier(425,135)(439,122)(439,104)
\qbezier(467,106)(467,87)(480,74)\qbezier(480,74)(493,61)(512,61)\qbezier(512,61)(530,61)(543,74)
\qbezier(543,74)(557,87)(557,106)\qbezier(467,106)(467,124)(480,137)\qbezier(480,137)(493,151)(512,151)
\qbezier(512,151)(530,151)(543,137)\qbezier(543,137)(557,124)(557,106)
\put(357,352){\circle*{4}}
\put(415,380){\circle*{4}}
\qbezier(414,391)(414,372)(426,359)\qbezier(426,359)(439,347)(458,347)\qbezier(458,347)(476,347)(489,359)
\qbezier(489,359)(502,372)(502,391)\qbezier(414,391)(414,409)(426,422)\qbezier(426,422)(439,435)(458,435)
\qbezier(458,435)(476,435)(489,422)\qbezier(489,422)(502,409)(502,391)
\put(2,237){\circle*{4}}
\put(195,216){\circle*{4}}
\qbezier(117,247)(117,228)(130,215)\qbezier(130,215)(143,202)(162,202)\qbezier(162,202)(180,202)(193,215)
\qbezier(193,215)(207,228)(207,247)\qbezier(117,247)(117,265)(130,278)\qbezier(130,278)(143,292)(162,292)
\qbezier(162,292)(180,292)(193,278)\qbezier(193,278)(207,265)(207,247)
\put(308,215){\circle*{4}}
\put(276,200){\circle*{4}}
\put(400,202){\circle*{4}}
\put(352,249){\circle*{4}}
\put(470,246){\circle*{4}}
\put(480,218){\circle*{4}}
\put(5,83){\circle*{4}}
\put(88,107){\circle*{4}}
\put(119,87){\circle*{4}}
\put(200,88){\circle*{4}}
\put(247,68){\circle*{4}}
\put(297,67){\circle*{4}}
\put(312,83){\circle*{4}}
\put(350,105){\circle*{4}}
\put(418,67){\circle*{4}}
\put(434,84){\circle*{4}}
\put(468,112){\circle*{4}}
\put(471,88){\circle*{4}}
\put(553,89){\circle*{4}}
\put(90,322){$G_{33}$}
\put(206,322){$G_{34}$}
\put(329,322){$G_{35}$}
\put(444,322){$G_{36}$}
\put(31,178){$G_{37}$}
\put(145,178){$G_{38}$}
\put(262,178){$G_{39}$}
\put(385,178){$G_{40}$}
\put(501,178){$G_{41}$}
\put(30,35){$G_{42}$}
\put(145,35){$G_{43}$}
\put(262,35){$G_{44}$}
\put(381,35){$G_{45}$}
\put(500,35){$G_{46}$}
\qbezier(62,399)(83,417)(104,434)
\put(102,347){\circle*{4}}
\qbezier(64,374)(105,373)(146,372)
\qbezier(64,374)(84,404)(104,434)
\put(223,433){\circle*{4}}
\qbezier(179,407)(201,420)(223,433)
\qbezier(223,433)(242,418)(262,403)
\qbezier(219,344)(221,389)(223,433)
\qbezier(219,344)(199,376)(179,407)
\qbezier(219,344)(240,374)(262,403)
\put(339,435){\circle*{4}}
\qbezier(296,406)(324,413)(339,435)
\qbezier(379,408)(352,418)(339,435)
\qbezier(295,378)(317,407)(339,435)
\qbezier(295,378)(337,393)(379,408)
\qbezier(295,378)(338,380)(381,381)
\qbezier(295,378)(326,365)(357,352)
\put(457,435){\circle*{4}}
\qbezier(417,409)(439,413)(457,435)
\qbezier(498,410)(471,414)(457,435)
\qbezier(444,350)(430,380)(417,409)
\qbezier(444,350)(450,393)(457,435)
\qbezier(444,350)(472,369)(501,387)
\qbezier(444,350)(471,380)(498,410)
\put(45,290){\circle*{4}}
\qbezier(45,290)(30,270)(5,264)
\qbezier(45,290)(64,267)(85,265)
\put(49,203){\circle*{4}}
\qbezier(2,237)(23,264)(45,290)
\qbezier(2,237)(25,220)(49,203)
\qbezier(2,237)(41,228)(80,219)
\qbezier(2,237)(45,242)(89,246)
\qbezier(2,237)(43,251)(85,265)
\put(161,292){\circle*{4}}
\put(168,203){\circle*{4}}
\qbezier(161,292)(147,265)(119,263)
\qbezier(204,264)(175,270)(161,292)
\qbezier(137,209)(149,251)(161,292)
\qbezier(137,209)(170,237)(204,264)
\qbezier(137,209)(172,224)(207,239)
\qbezier(195,216)(166,213)(137,209)
\qbezier(137,209)(128,236)(119,263)
\put(278,289){\circle*{4}}
\qbezier(278,289)(262,268)(236,266)
\qbezier(316,265)(290,266)(278,289)
\qbezier(276,200)(277,245)(278,289)
\qbezier(276,200)(254,221)(232,241)
\qbezier(276,200)(256,233)(236,266)
\qbezier(276,200)(296,233)(316,265)
\qbezier(276,200)(298,221)(320,241)
\put(395,289){\circle*{4}}
\put(431,218){\circle*{4}}
\qbezier(395,289)(383,268)(357,267)
\qbezier(435,267)(411,261)(395,289)
\qbezier(395,289)(376,248)(352,249)
\qbezier(352,249)(391,234)(431,218)
\qbezier(352,249)(402,223)(400,202)
\qbezier(352,249)(380,224)(371,208)
\qbezier(352,249)(396,247)(440,245)
\qbezier(352,249)(393,258)(435,267)
\put(513,292){\circle*{4}}
\put(499,206){\circle*{4}}
\qbezier(513,292)(500,266)(475,269)
\qbezier(513,292)(524,267)(552,270)
\qbezier(480,218)(490,254)(475,269)
\qbezier(480,218)(514,221)(525,205)
\qbezier(480,218)(528,228)(550,222)
\qbezier(480,218)(519,234)(558,249)
\qbezier(480,218)(496,255)(513,292)
\qbezier(480,218)(516,244)(552,270)
\put(45,149){\circle*{4}}
\qbezier(0,105)(0,86)(12,73)\qbezier(12,73)(25,61)(44,61)\qbezier(44,61)(62,61)(75,73)\qbezier(75,73)(88,86)(88,105)
\qbezier(0,105)(0,123)(12,136)\qbezier(12,136)(25,149)(44,149)\qbezier(44,149)(62,149)(75,136)
\qbezier(75,136)(88,123)(88,105)
\qbezier(81,128)(54,122)(45,149)
\put(58,63){\circle*{4}}
\qbezier(26,63)(35,106)(45,149)
\qbezier(26,63)(53,96)(81,128)
\qbezier(26,63)(65,82)(82,82)
\qbezier(26,63)(57,85)(88,107)
\qbezier(45,149)(33,130)(6,129)
\qbezier(6,129)(27,114)(26,63)
\qbezier(0,107)(15,102)(26,63)
\put(162,149){\circle*{4}}
\qbezier(162,149)(148,125)(122,129)
\qbezier(162,149)(176,124)(197,129)
\put(187,71){\circle*{4}}
\put(135,68){\circle*{4}}
\qbezier(116,108)(149,115)(162,149)
\qbezier(116,108)(133,89)(135,68)
\qbezier(116,108)(139,84)(166,62)
\qbezier(116,108)(151,90)(187,71)
\qbezier(116,108)(158,98)(200,88)
\qbezier(116,108)(159,109)(203,110)
\qbezier(116,108)(156,119)(197,129)
\put(278,147){\circle*{4}}
\put(268,61){\circle*{4}}
\qbezier(278,147)(264,128)(235,127)
\qbezier(314,124)(291,123)(278,147)
\qbezier(233,85)(252,113)(235,127)
\qbezier(233,85)(282,79)(297,67)
\qbezier(233,85)(286,86)(312,83)
\qbezier(233,85)(275,95)(317,104)
\qbezier(233,85)(273,105)(314,124)
\qbezier(233,85)(255,116)(278,147)
\qbezier(268,61)(269,72)(233,85)
\put(394,148){\circle*{4}}
\put(394,60){\circle*{4}}
\qbezier(394,148)(377,126)(355,126)
\qbezier(394,148)(410,126)(432,127)
\qbezier(369,67)(381,108)(394,148)
\qbezier(369,67)(400,97)(432,127)
\qbezier(369,67)(403,86)(438,105)
\qbezier(369,67)(363,96)(350,105)
\qbezier(369,67)(370,110)(355,126)
\qbezier(369,67)(405,76)(418,67)
\qbezier(434,84)(409,82)(369,67)
\put(515,150){\circle*{4}}
\qbezier(515,150)(495,128)(476,133)
\qbezier(515,150)(528,131)(547,133)
\put(511,61){\circle*{4}}
\qbezier(511,61)(513,106)(515,150)
\qbezier(511,61)(493,97)(476,133)
\qbezier(511,61)(529,97)(547,133)
\qbezier(511,61)(489,87)(468,112)
\qbezier(511,61)(533,88)(556,114)
\qbezier(511,61)(486,85)(471,88)
\qbezier(511,61)(537,83)(553,89)
\put(205,-9){Fig. 3.5. $G_{33}$-$G_{46}$}
\end{picture}
\end{center}

\bigskip

\begin{center}
 \begin{tabular}{|c|c|c|c|c|c|c|c|c|}
\hline graph & $Q$-index & graph & $Q$-index & graph & $Q$-index & graph & $Q$-index \\
\hline $G_{33}$  & $7.4035$  & $G_{37}$   &  $9.0193$ & $G_{41}$ & $9.8476$       & $G_{45}$  & $10.7002$  \\
\hline $G_{34}$  & $7.8845$  & $G_{38}$   &  $9.0704$ & $G_{42}$ & $9.8521$       & $G_{46}$     & $10.7005$ \\
\hline $G_{35}$  & $8.3281$  & $G_{39}$   &  $7.4621$ & $G_{43}$ & $10.6779$       &    &  \\
\hline $G_{36}$  & $8.4076$  & $G_{40}$   &  $9.8162$ & $G_{44}$ & $10.6976$        &     &    \\
 \hline

\end{tabular}
\end{center}
\begin{center}
{Table 3. The approximation of the $Q$-index for $G_{i}$ $(33\leq i\leq 46)$}
\end{center}

When $n=8$, $G$ is isomorphic to one in $\{G_{33}, G_{34}\}$; when $n=9$, $G$ is isomorphic to one in $\{G_{35}, G_{36}\}$; when $n=10$, $G$ is isomorphic to one in $\{G_{37}, G_{38}, G_{39}\}$; when $n=11$, $G$ is isomorphic to one in $\{G_{40}, G_{41}, G_{42}\}$; when $n=12$, $G$ is isomorphic to one in $\{G_{43}, G_{44}, G_{45}, G_{46}\}$ (see Fig. 3.5). By computation with computer, we get the approximation of the $Q$-index for each $G_{i}$ $(33\leq i\leq 46)$ (see Table 3). And by computation with computer, we get that when $n=12$, $q(\mathcal {H})\approx 12.4233$. Combining with the results about $q(\mathcal {H})$ for order $n=8$, $9$, $10$, $11$ in Case 1 and in the proof of Lemma \ref{le3,2}, by a simple comparison, we get that for each $G_{i}$ ($33\leq i\leq 46$) of order $n$ ($8\leq n\leq 12$), $q(G_{i})<q(\mathcal {H})$.

{\bf Case 3} $G\cong D_{3}$.

For this case, $n\geq 7$. For convenience, we suppose that $G= D_{3}$. Because $\Delta(G) = n -3$, $k\notin \{0$, $1$, $j+1$, $j+2$, $j+3$, $n-1\}$. Let $F=G-v_{1}v_{n-1}+v_{k}v_{0}$.  As Lemma \ref{le3,2}, it can be proved that if $n\geq 10$, then $q(G)<q(F)$. By Lemma \ref{le3,2}, we get that $q(F)<q(\mathcal {H})$. Then $q(G)<q(\mathcal {H})$.

\setlength{\unitlength}{0.63pt}
\begin{center}
\begin{picture}(458,294)
\qbezier(1,244)(1,224)(15,210)\qbezier(15,210)(29,196)(49,196)\qbezier(49,196)(68,196)(82,210)
\qbezier(82,210)(97,224)(97,244)\qbezier(1,244)(1,263)(15,277)\qbezier(15,277)(29,292)(49,292)
\qbezier(49,292)(68,292)(82,277)\qbezier(82,277)(97,263)(97,244)
\qbezier(188,245)(188,225)(202,211)\qbezier(202,211)(216,197)(236,197)\qbezier(236,197)(255,197)(269,211)
\qbezier(269,211)(284,225)(284,245)\qbezier(188,245)(188,264)(202,278)\qbezier(202,278)(216,293)(236,293)
\qbezier(236,293)(255,293)(269,278)\qbezier(269,278)(284,264)(284,245)
\qbezier(362,246)(362,226)(376,212)\qbezier(376,212)(390,198)(410,198)\qbezier(410,198)(429,198)(443,212)
\qbezier(443,212)(458,226)(458,246)\qbezier(362,246)(362,265)(376,279)\qbezier(376,279)(390,294)(410,294)
\qbezier(410,294)(429,294)(443,279)\qbezier(443,279)(458,265)(458,246)
\qbezier(0,92)(0,71)(14,57)\qbezier(14,57)(28,43)(49,43)\qbezier(49,43)(69,43)(83,57)\qbezier(83,57)(98,71)(98,92)
\qbezier(0,92)(0,112)(14,126)\qbezier(14,126)(28,141)(49,141)\qbezier(49,141)(69,141)(83,126)
\qbezier(83,126)(98,112)(98,92)
\qbezier(188,93)(188,73)(202,59)\qbezier(202,59)(216,45)(236,45)\qbezier(236,45)(255,45)(269,59)
\qbezier(269,59)(284,73)(284,93)\qbezier(188,93)(188,112)(202,126)\qbezier(202,126)(216,141)(236,141)
\qbezier(236,141)(255,141)(269,126)\qbezier(269,126)(284,112)(284,93)
\qbezier(362,94)(362,74)(376,60)\qbezier(376,60)(390,46)(410,46)\qbezier(410,46)(429,46)(443,60)
\qbezier(443,60)(458,74)(458,94)\qbezier(362,94)(362,113)(376,127)\qbezier(376,127)(390,142)(410,142)
\qbezier(410,142)(429,142)(443,127)\qbezier(443,127)(458,113)(458,94)
\put(31,288){\circle*{4}}
\put(1,237){\circle*{4}}
\qbezier(31,288)(16,263)(1,237)
\put(63,289){\circle*{4}}
\put(97,241){\circle*{4}}
\qbezier(63,289)(80,265)(97,241)
\put(88,271){\circle*{4}}
\put(8,269){\circle*{4}}
\put(18,208){\circle*{4}}
\qbezier(18,208)(24,248)(31,288)
\qbezier(18,208)(57,225)(97,241)
\put(71,202){\circle*{4}}
\qbezier(18,208)(40,249)(63,289)
\put(201,276){\circle*{4}}
\put(201,213){\circle*{4}}
\qbezier(201,276)(201,245)(201,213)
\put(271,276){\circle*{4}}
\put(271,214){\circle*{4}}
\qbezier(271,276)(271,245)(271,214)
\put(237,292){\circle*{4}}
\put(237,198){\circle*{4}}
\qbezier(237,198)(219,237)(201,276)
\qbezier(237,198)(237,245)(237,292)
\qbezier(237,198)(254,237)(271,276)
\put(188,245){\circle*{4}}
\put(284,246){\circle*{4}}
\put(389,289){\circle*{4}}
\put(363,236){\circle*{4}}
\qbezier(389,289)(376,263)(363,236)
\put(431,288){\circle*{4}}
\put(456,243){\circle*{4}}
\qbezier(431,288)(443,266)(456,243)
\put(442,212){\circle*{4}}
\put(412,198){\circle*{4}}
\put(379,211){\circle*{4}}
\qbezier(379,211)(384,250)(389,289)
\qbezier(379,211)(405,250)(431,288)
\qbezier(379,211)(417,227)(456,243)
\qbezier(379,211)(410,212)(442,212)
\put(366,266){\circle*{4}}
\put(453,268){\circle*{4}}
\put(26,134){\circle*{4}}
\put(1,87){\circle*{4}}
\qbezier(26,134)(13,111)(1,87)
\put(70,136){\circle*{4}}
\put(98,87){\circle*{4}}
\qbezier(70,136)(84,112)(98,87)
\put(14,58){\circle*{4}}
\put(49,44){\circle*{4}}
\put(80,55){\circle*{4}}
\qbezier(49,44)(37,89)(26,134)
\qbezier(49,44)(59,90)(70,136)
\qbezier(49,44)(73,66)(98,87)
\qbezier(49,44)(25,66)(1,87)
\put(205,129){\circle*{4}}
\put(190,80){\circle*{4}}
\qbezier(205,129)(197,105)(190,80)
\put(267,128){\circle*{4}}
\put(280,79){\circle*{4}}
\qbezier(267,128)(273,104)(280,79)
\put(238,141){\circle*{4}}
\put(394,49){\circle*{4}}
\put(258,51){\circle*{4}}
\put(219,49){\circle*{4}}
\qbezier(219,49)(212,89)(205,129)
\qbezier(219,49)(228,95)(238,141)
\qbezier(219,49)(249,64)(280,79)
\put(376,127){\circle*{4}}
\put(367,75){\circle*{4}}
\put(443,128){\circle*{4}}
\put(454,77){\circle*{4}}
\qbezier(443,128)(448,103)(454,77)
\put(434,53){\circle*{4}}
\put(413,141){\circle*{4}}
\qbezier(413,141)(403,95)(394,49)
\qbezier(413,141)(423,97)(434,53)
\qbezier(413,141)(390,108)(367,75)
\qbezier(413,141)(433,109)(454,77)
\put(457,103){\circle*{4}}
\put(363,103){\circle*{4}}
\put(282,105){\circle*{4}}
\put(190,109){\circle*{4}}
\put(92,116){\circle*{4}}
\put(6,117){\circle*{4}}
\put(33,172){$G_{47}$}
\put(221,171){$G_{48}$}
\put(398,171){$G_{49}$}
\put(36,20){$G_{50}$}
\put(224,21){$G_{51}$}
\put(399,22){$G_{52}$}
\put(160,-9){Fig. 3.6. $G_{47}$-$G_{52}$}
\qbezier(376,127)(371,101)(367,75)
\qbezier(219,49)(243,89)(267,128)
\end{picture}
\end{center}

\begin{center}
 \begin{tabular}{|c|c|c|c|c|c|c|}
\hline graph & $Q$-index & graph & $Q$-index & graph & $Q$-index \\
\hline $G_{47}$  & $7.6044$  & $G_{49}$   &  $8.2339$ & $G_{51}$ & $8.2078$    \\
\hline $G_{48}$  & $7.4741$  & $G_{50}$   &  $8.2833$ & $G_{52}$ & $8.1408$    \\ \hline

\end{tabular}
\end{center}
\begin{center}
{Table 4. The approximation of the $Q$-index for $G_{i}$ $(1\leq i\leq 9)$}
\end{center}

When $n=7$, then $G\cong G_{10}$. From Case 1, we know that for $n=7$, $q(G_{10})<q(\mathcal {H})$.  when $n=8$, $G$ is isomorphic to one in $\{G_{47}, G_{48}\}$; when $n=9$, $G$ is isomorphic to one in $\{G_{49}, G_{50}, G_{51}, G_{52}\}$ (see Fig. 3.6). By computation with computer, we get the approximation of the $Q$-index for each $G_{i}$ $(47\leq i\leq 52)$ (see Table 4). Combining with the results about $q(\mathcal {H})$ for order $n=8$, $9$ in the proof of Lemma \ref{le3,2}, by a simple comparison, we get that for each $G_{i}$ ($i=10$ and $47\leq i\leq 52$) of order $n$ ($n= 8, 9$), $q(G_{i})<q(\mathcal {H})$.

From above three cases, it follows that for a maximal outer-planar graph $G$ with order $n \geq 7$ and $\Delta(G) = n -3$, $q(G)<q(\mathcal {H})$. This completes the proof.
\ \ \ \ \ $\Box$
\end{proof}

\begin{lemma}\label{le3,4} 
Let $G$ a maximal outer-planar graph of order $n$. Then $q(G)\leq q(\mathcal {H})$ with equality if and only if $G\cong \mathcal {H}$.
\end{lemma}

\begin{proof}
It can be checked that when $n=1$, $2$, $3$, $4$, $5$, $G\cong \mathcal {H}$.
Then the result follows from Lemmas \ref{le3,1}-\ref{le3,3}. This completes the proof.
\ \ \ \ \ $\Box$
\end{proof}

\begin{theorem}\label{th3,5} 
Let $G$ be an outer-planar graph of order $n$. Then $q(G)<q(\mathcal {H})$ with equality if and only if $G\cong \mathcal {H}$.
\end{theorem}

\begin{proof}
From the narration in Section 2, we know that by adding edges, a maximal outer-planar graph can be obtained from a nonmaximal outer-planar graph; and know that for a graph $G$, $q(G+e)> q(G)$ if $e\notin E(G)$.
Then the result follows from Lemma \ref{le3,4}. This completes the proof.
\ \ \ \ \ $\Box$
\end{proof}

\small {

}

\end{document}